\documentclass[12pt]{article}

\usepackage[utf8]{inputenc}
\usepackage[affil-it]{authblk}
\usepackage{amsfonts, amsmath, amssymb, amsthm, bbm, color, enumerate, graphicx, mathtools, tikz, hyperref, relsize, bm, adjustbox, multirow}
\usepackage[style=authoryear]{biblatex}
\usepackage{caption}
\usepackage[margin=30mm]{geometry}

\newtheorem{theorem}{\bf Theorem}

\newtheorem{lemma}{Lemma}
\newtheorem{redlemma}[lemma]{\textcolor{black}{Lemma}}

\newtheorem{remark}{Remark}

\DeclareMathOperator*{\argmax}{arg\,max}

\newcommand{\Var}{\mathrm{Var}}

\allowdisplaybreaks

\title{Asymptotic Theory for Estimation of the H\"usler-Reiss Distribution via Block Maxima Method}
\author[1]{Hank Flury}
\author[1]{Jan Hannig}
\author[1]{Richard Smith}
\affil[1]{University of North Carolina at Chapel Hill, Department of Statistics and Operations Research, Chapel Hill, North Carolina, United States}

\begin{document}
\setcounter{lemma}{1}
\maketitle
Hank Flury is the corresponding author. His email is: fluryh@unc.edu. Hank Flury's ORCiD is 0009-0004-1722-8607. Richard Smith's ORCiD is 0000-0002-3454-3199. Jan Hannig's ORCiD is 0000-0002-4164-0173.

Jan Hannig's research was supported in part by the National Science Foundation under Grant No. DMS-1916115, 2113404, and 2210337. Hank Flury was also supported in part by the National Science Foundation under Grant No. DMS-2210337.
\newpage

\begin{abstract}
    The H\"usler-Reiss distribution describes the limit of the pointwise maxima of a bivariate normal distribution. This distribution is defined by a single parameter, $\lambda$. We provide asymptotic theory for maximum likelihood estimation of $\lambda$ under a block maxima approach. Our work assumes independent and identically distributed bivariate normal random variables, grouped into blocks where the block size and number of blocks increase simultaneously. With these assumptions our results provide conditions for the asymptotic normality of the Maximum Likelihood Estimator (MLE). We characterize the bias of the MLE, provide conditions under which this bias is asymptotically negligible, and discuss how to choose the block size to minimize a bias-variance trade-off. The proofs are an extension of previous results for choosing the block size in the estimation of univariate extreme value distributions (\cite{dombry19}), providing a potential basis for extensions to multivariate cases where both the marginal and dependence parameters are unknown. The proofs rely on the Argmax Theorem applied to a localized log-likelihood function, combined with a Lindeberg-Feller Central Limit Theorem argument to establish asymptotic normality. Possible applications of the method include composite likelihood estimation in Brown-Resnick processes, where it is known that the bivariate distributions are of H\"usler-Reiss form. 
\end{abstract}

\section{Introduction}
The Block Maxima (BM) method is one of the main techniques in Extreme Value Theory (EVT), with the goal of estimation of the limiting distribution of the maxima of a random variable. In this context, the distribution is limiting in the sense that we take the maxima over an increasingly large number of samples. In this method, we take $n$ independent and identically distributed (i.i.d.) observations and split them into $k_n$ equally sized groups of $m_n$ samples. Going forward, we drop the subscript $n$ for the sake of readability. Then, we can take the maximum over our $m$ data points for each block, now having $k$ values with which to do inference. Two questions naturally arise from this practice: If we are trying to estimate a limiting distribution, and our sample comes from a distribution taken under a finite sample size, is it reasonable to estimate the parameters of the limiting distribution with our sample? Assuming that it is reasonable, how should a practitioner of the BM method split their $n$ samples into $k$ and $m$?

\textcites{dombry19} answers these questions for the univariate case. They showed that, under mild assumptions on the distribution of the sample and the relationship between $k$ and $m$, with probability tending toward 1, there exists a maximum likelihood estimator, $\hat{\theta}$, such that we have a Central Limit Theorem-like result:
\begin{equation}
    \sqrt{k}\left(\hat{\theta} - \theta\right) \xrightarrow{d} N\left(I_{\theta}^{-1}\lambda, I_{\theta}^{-1}\right)
\end{equation}
where $I_{\theta}$ is the Fisher Information of the limiting distribution with true parameters $\theta$, and $\lambda$ is a parameter that arises as the limit of the relationship between $k$ and $m$. This work will develop similar results for a specific bivariate case. One of the complications of the bivariate case is that we do not have a general form for the limiting distribution, while in the univariate case we know that the limiting distribution is a generalized extreme value (GEV) distribution. Thus, for our work, we need to assume a limiting distribution.

The H\"usler-Reiss distribution arises as the limit of properly scaled maxima of bivariate normal data (\cite{huslerreiss89}). The distribution is given by
\begin{equation}
    H_\lambda\left(x,y\right) = \exp\left\{-e^{-x}\Phi\left(\lambda + \frac{y-x}{2\lambda}\right)-e^{-y}\Phi\left(\lambda + \frac{x-y}{2\lambda}\right)\right\}.
\end{equation}
where $\Phi\left(t\right)$ is the distribution function of a standard normal random variable. Namely, if $\left(X_i, Y_i\right)$ are i.i.d. bivariate standard normal random variables with correlation $\rho_m$, and we define $b_m$ to satisfy the equation $b_m = m \phi\left(b_m\right)$, where $\phi(x) = \frac{1}{\sqrt{2\pi}}e^{-\frac{1}{2}x^2}$ is the standard normal density, then the H\"usler-Reiss distribution will arise as the distribution of $\left(b_m\max\left(X_1,...,X_m\right)-b_m^2,b_m\max\left(Y_1,...Y_m\right)-b_m^2\right)$ as $m \rightarrow \infty$. The parameter $\lambda \in \left[0,\infty\right]$ arises as the limit of $\left(1-\rho_{m}\right)\ln\left(m\right) \rightarrow \lambda^2$, with the convention that $\lambda = 0$ corresponds to complete dependence (making the H\"usler-Reiss distribution equivalent to a univariate Gumbel distribution with $X_i = Y_i$) and $\lambda = \infty$ corresponds to complete independence (making the H\"usler-Reiss distribution equivalent to the product of two independent Gumbel distributions). \textcolor{black}{The bivariate H\"usler-Reiss distribution is often used to model extreme climactic events, where the two variables represent the magnitude of such an event at two sites. The worst effects of such extreme climactic events are a result of simultaneous occurrences of extreme and the likelihood of such concurrent extremes hinges on dependence parameters like $\lambda$. As such, any conditional predictions, say the likelihood of extreme rainfall in one area given extreme rainfall at a nearby site requires accurate estimation of the dependence parameter. Indeed the issue of fitting a H\"usler-Reiss distribution is well documented in extreme value literature. \textcite{davisonpadoanribatet} consider the H\"usler-Reiss model for extreme rainfall, and \textcite{Asadi_2015} for river flooding. Any study in which estimation of such distributions is desired will hinge on accurate estimation of the $\lambda$ parameter.}

While the H\"usler-Reiss distribution is useful in context other than bivariate maxima, this is the only case that is studied in the literature. Our work continues to focus on this framework for possible extension to other cases in the future. Maximum likelihood estimation of the H\"usler-Reiss distribution was discussed by \textcites{engelke2015}. Their work, however, conditions their estimator on the observation of an extreme event, leading to entirely different results from those given here. In contrast, \textcites{tawn88} and \textcites{tawn90} derive results for the MLE of dependence parameters of bivariate extreme values (not including the H\"usler-Reiss distribution), focusing on properties of the MLE when the true value is on the boundaries of the parameter space, corresponding to independence or exact dependence. Our work focuses on the properties of the MLE when $\lambda \in (0,\infty)$, leading to new results. The uniform convergence rate of the maxima of scaled bivariate normal random variables was derived by \textcites{liao}. Our work extends the probabilistic results of \textcite{liao} to the statistical estimation problem. 

Our work is broken into the following sections. Section 2 discusses the relevant details of maximum likelihood estimation in the case that our data are distributed exactly according to a H\"usler-Reiss distribution. Section 3 handles the case of the misspecified model, where our data are the scaled maxima of bivariate normal random variables but we approximate this joint distribution by the H\"usler-Reiss distribution and try to estimate $\lambda$. Section 4 provides a simulation study to support our theoretical findings. Finally, the appendix contains details of the proofs.

\section{Estimation Under the H\"usler-Reiss Model} \label{Estimation Under the HR Model}
We begin with some properties of the H\"usler-Reiss model, the main result being a rigorous proof that the Fisher Information is finite. Although we do not prove the existence of asymptotic normality of the MLE for $\lambda$ when our data are distributed according to a H\"usler-Reiss distribution as a separate case, we note that this arises as a consequence of Theorem~\ref{MLE Form}. It will be useful later to have analyzed the log-likelihood and know its form. We will denote $G\left(x,y,\lambda\right) := e^{-x}\Phi\left(\lambda +\frac{y-x}{2\lambda}\right) + e^{-y}\Phi\left(\lambda +\frac{x-y}{2\lambda}\right)$ so that $H_{\lambda}\left(x,y\right) = e^{-G\left(x,y,\lambda\right)}$. Where necessary for space, we omit the arguments of $G(x,y,\lambda)$ when there is no risk of ambiguity. We begin by finding the form of the derivative of the log-likelihood,
\begin{align*}
    &l'\left(x,y,\lambda\right) := \sum_{i=1}^k\frac{\partial}{\partial\lambda}\ln\left(\frac{\partial^2}{\partial x\partial y}H_\lambda\left(x_i,y_i\right)\right)\\ 
    &= \sum_{i=1}^k\frac{\partial}{\partial\lambda}\ln\left(e^{-G\left(x_i,y_i,\lambda\right)}\left(\frac{\partial}{\partial y}G\left(x_i,y_i,\lambda\right)\frac{\partial}{\partial x}G\left(x_i,y_i,\lambda\right) -\frac{\partial^2}{\partial y\partial x}G\left(x_i,y_i,\lambda\right)\right)\right)\\
    &=\sum_{i=1}^k\frac{\frac{\partial}{\partial \lambda} \left[\frac{\partial}{\partial y}G\frac{\partial}{\partial x}G\right] - \frac{\partial^3}{\partial y\partial x\partial\lambda}G - \left(\frac{\partial}{\partial\lambda}G\right)\left(\frac{\partial}{\partial y}G\frac{\partial}{\partial x}G -\frac{\partial^2}{\partial y\partial x}G\right)}{\frac{\partial}{\partial y}G\frac{\partial}{\partial x}G -\frac{\partial^2}{\partial y\partial x}G}.
\end{align*}

We calculate the forms of each term to be:

\begin{align*}
    \frac{\partial}{\partial\lambda}G\left(x,y,\lambda\right) &= 2e^{-x}\phi\left(\lambda + \frac{y-x}{2\lambda}\right),\\
    \frac{\partial}{\partial x}G\left(x,y,\lambda\right) &= -e^{-x}\Phi\left(\lambda  +\frac{y-x}{2\lambda}\right),\\
    \frac{\partial}{\partial y}G\left(x,y,\lambda\right) &= -e^{-y}\Phi\left(\lambda  +\frac{x-y}{2\lambda}\right),\\
     \frac{\partial^2}{\partial x \partial y}G\left(x,y,\lambda\right) &= -e^{-x}\phi\left(\lambda + \frac{y-x}{2\lambda}\right)\frac{1}{2\lambda},\\
     \frac{\partial^3}{\partial x\partial y\partial\lambda}G\left(x,y,\lambda\right) &= e^{-x}\phi\left(\lambda+\frac{y-x}{2\lambda}\right)\left(\frac{1}{2\lambda^2}+\frac{1}{2}-\frac{\left(y-x\right)^2}{8\lambda^4}\right),\\
     \frac{\partial}{\partial\lambda}\left[\frac{\partial}{\partial x}G\left(x,y,\lambda\right)\frac{\partial}{\partial y}G\left(x,y,\lambda\right)\right]
     &=e^{-x}\phi\left(\lambda + \frac{y-x}{2\lambda}\right)\left[\Phi\left(\lambda+\frac{x-y}{2\lambda}\right)e^{-y}\left(1-\frac{y-x}{2\lambda^2}\right)\right.\\
     &\hspace{10mm}+ \left.\Phi\left(\lambda+\frac{y-x}{2\lambda}\right)e^{-x}\left(1-\frac{x-y}{2\lambda^2}\right)\right].
\end{align*}
Our first goal is to show that the Fisher Information for this distribution exists. To do this we begin with a lemma:
\begin{lemma}\label{integral}
\begin{equation}
    \int_{-\infty}^{\infty}\int_{-\infty}^{\infty}e^{-ax}e^{-by}\phi\left(\lambda + \frac{y-x}{2\lambda}\right)x^{k_1}y^{k_2}\left(x-y\right)^{k_3}e^{-e^{-x}f\left(x-y\right)}dxdy
\end{equation}
is finite for $a,b\in \mathbb{Z}$, such that $a+b \geq 1,  \lambda > 0, k_i \in \mathbb{N} \cup \left\{0\right\}, i =1,2,3$ and $f$ such that there exists $c>0$ such that $c < f\left(t\right)$ for all real $t$.
\end{lemma}

To show that this integral is finite, we show that the integral of the absolute value of the integrand is finite. It can be shown through integration by parts and induction that $\int_{-\infty}^{\infty}e^{-kx}e^{-ce^{-x}}dx = \frac{\left(k-1\right)!}{c^k}$. Using the change of variable $\delta = x-y$,
\begin{align*}
    0 &\leq \int_{-\infty}^{\infty}\int_{-\infty}^{\infty}\left|e^{-ax}e^{-by}\phi\left(\lambda + \frac{y-x}{2\lambda}\right)x^{k_1}y^{k_2}\left(x-y\right)^{k_3}e^{-e^{-x}f\left(x-y\right)}\right|dxdy \\
    &\leq \sum_{i=0}^{k_2}\binom{k_2}{i}\int_{-\infty}^{\infty}\left|\delta^{k_3+k_2-i}\right|e^{b\delta}\phi\left(\lambda - \frac{\delta}{2\lambda}\right)d\delta\int_{-\infty}^{\infty}\left|x^{k_1+i}\right|e^{-\left(a+b\right)x}e^{-ce^{-x}}dx.
\end{align*}
Note that the first integral depends only on the integer moments of a normal distribution, which must be finite. To see that the second integral is finite, we use the bounds $|x^t| < e^{-2tx}$ when $x \leq 0$. This gives us

\begin{align*}
    \int_{-\infty}^{\infty}|x^{k_1+i}|e^{-\left(a+b\right)x}e^{-ce^{-x}}dx &\leq \int_{-\infty}^0e^{-\left(2k_1+2i+a+b\right)x}e^{-ce^{-x}}dx + \int_{0}^{\infty}x^{k_1+i}e^{-\left(a+b\right)x}dx\\
    &\leq \frac{\left(2k_1+2i+a+b\right)!}{c^{2k_1+2i+a+b}} + \left(a+b\right)^{-1-k_1-i}\Gamma\left(1 + k_1 + i\right),
\end{align*}
which is surely finite.
\qed

With this lemma in hand, we are ready to prove that the Fisher Information must exist. 

\begin{theorem}\label{MLE existence}
The Fisher Information for a H\"usler-Reiss distribution with parameter $\lambda \in \left(0,\infty\right)$ is finite.
\end{theorem}
\noindent \textit{Proof.}

The proof of Theorem~\ref{MLE existence} involves showing that each term of the integrand corresponding to the Fisher Information is bounded by a form of the integrand above, up to a constant. We take a derivative of $l'\left(x,y,\lambda\right)$ with respect to $\lambda$ to find:

\begin{align}\label{second derivative form}
    &\frac{\partial}{\partial\lambda}l'\left(x,y,\lambda\right) = \frac{\partial^2}{\partial^2\lambda}\ln\left(\frac{\partial^2}{\partial x\partial y}H_{\lambda}\left(x,y\right)\right)\nonumber\\
    =&-\frac{\partial^2}{\partial^2\lambda}G\left(x,y,\lambda\right)+ \frac{\frac{\partial^2}{\partial^2 \lambda} \left[\frac{\partial}{\partial y}G\left(x,y,\lambda\right)\frac{\partial}{\partial x}G\left(x,y,\lambda\right)\right] - \frac{\partial^4}{\partial y\partial x\partial^2\lambda}G\left(x,y,\lambda\right)}{\frac{\partial}{\partial y}G\left(x,y,\lambda\right)\frac{\partial}{\partial x}G\left(x,y,\lambda\right) -\frac{\partial^2}{\partial y\partial x}G\left(x,y,\lambda\right)}\nonumber\\
    &- \left(\frac{\frac{\partial}{\partial \lambda} \left[\frac{\partial}{\partial y}G\left(x,y,\lambda\right)\frac{\partial}{\partial x}G\left(x,y,\lambda\right)\right] - \frac{\partial^3}{\partial y\partial x\partial\lambda}G\left(x,y,\lambda\right)}{\frac{\partial}{\partial y}G\left(x,y,\lambda\right)\frac{\partial}{\partial x}G\left(x,y,\lambda\right) -\frac{\partial^2}{\partial y\partial x}G\left(x,y,\lambda\right)}\right)^2.
\end{align}

We begin with $-\frac{\partial^2}{\partial^2\lambda}G\left(x,y,\lambda\right) = 2e^{-x}\phi\left(\lambda+\frac{y-x}{2\lambda}\right)\left(\lambda-\frac{\left(y-x\right)^2}{4\lambda^3}\right)$. We see that our integrand is:
\begin{align*}
    &\int_{-\infty}^{\infty}\int_{-\infty}^{\infty}2e^{-x}\phi\left(\lambda+\frac{y-x}{2\lambda}\right)\left(\lambda-\frac{\left(y-x\right)^2}{4\lambda^3}\right)e^{{-e^{-x}\Phi\left(\lambda + \frac{y-x}{2\lambda}\right)-e^{-y}\Phi\left(\lambda + \frac{x-y}{2\lambda}\right)}}\\
    &\hspace{10mm} \left[e^{-x-y}\Phi\left(\lambda + \frac{y-x}{2\lambda}\right)\Phi\left(\lambda + \frac{x-y}{2\lambda}\right) + e^{-x}\phi\left(\lambda+\frac{y-x}{2\lambda}\right)\frac{1}{2\lambda}\right]dxdy\\
    &=\int_{-\infty}^{\infty}\int_{-\infty}^{\infty}2e^{-2x}\phi\left(\lambda-\frac{\delta}{2\lambda}\right)\left(\lambda-\frac{\delta^2}{4\lambda^3}\right)e^{{-e^{-x}\left(\Phi\left(\lambda - \frac{\delta}{2\lambda}\right)+e^{\delta}\Phi\left(\lambda + \frac{\delta}{2\lambda}\right)\right)}}\\
    &\hspace{10mm} \left[e^{\delta-x}\Phi\left(\lambda - \frac{\delta}{2\lambda}\right)\Phi\left(\lambda + \frac{\delta}{2\lambda}\right) + \phi\left(\lambda-\frac{\delta}{2\lambda}\right)\frac{1}{2\lambda}\right]dxd\delta.
\end{align*}

When factored, both terms of this integral are bounded by a form of the integral in Lemma~\ref{integral} with appropriate $a$, $b$, and $k_i, i=1,2,3$, since $0<\Phi\left(t\right),\phi\left(t\right)<1$ for all real $t$ and $\frac{1}{2} < \Phi\left(\lambda - \frac{\delta}{2\lambda}\right)+e^{\delta}\Phi\left(\lambda + \frac{\delta}{2\lambda}\right)$ for all values of $\delta$ since $\lambda > 0$. Thus we are assured that the expectation of this term is finite. We now move to the term $\frac{\frac{\partial^2}{\partial^2 \lambda} \left[\frac{\partial}{\partial y}G\left(x,y,\lambda\right)\frac{\partial}{\partial x}G\left(x,y,\lambda\right)\right]}{\frac{\partial}{\partial y}G\left(x,y,\lambda\right)\frac{\partial}{\partial x}G\left(x,y,\lambda\right) -\frac{\partial^2}{\partial y\partial x}G\left(x,y,\lambda\right)}$. We begin by finding the form of $\frac{\partial^2}{\partial^2 \lambda} \left[\frac{\partial}{\partial y}G\left(x,y,\lambda\right)\frac{\partial}{\partial x}G\left(x,y,\lambda\right)\right]$:
\begin{align*}
    &\frac{\partial^2}{\partial^2 \lambda} \left[\frac{\partial}{\partial y}G\left(x,y,\lambda\right)\frac{\partial}{\partial x}G\left(x,y,\lambda\right)\right]\\
    &= e^{-x}\phi\left(\lambda + \frac{y-x}{2\lambda}\right)\left[2e^{-x}\phi\left(\lambda + \frac{y-x}{2\lambda}\right)\left(1-\frac{\left(x-y\right)^2}{2\lambda^2}\right)\right.\\
    &\hspace{10mm}+ e^{-x}\Phi\left(\lambda + \frac{y-x}{2\lambda}\right)\left\{-\lambda + \frac{x-y}{2\lambda}+\frac{\left(y-x\right)^2}{4\lambda^3}-\frac{\left(x-y\right)^3}{8\lambda^5}+\frac{x-y}{\lambda^3}\right\}\\
    &\hspace{10mm}\left.+e^{-y}\Phi\left(\lambda + \frac{x-y}{2\lambda}\right)\left\{-\lambda + \frac{y-x}{2\lambda}+\frac{\left(x-y\right)^2}{4\lambda^3}-\frac{\left(y-x\right)^3}{8\lambda^5}+\frac{y-x}{\lambda^3}\right\}\right].
\end{align*}

We note that our density takes the form $e^{-G\left(x,y\right)}\left[\frac{\partial}{\partial y}G\left(x,y,\lambda\right)\frac{\partial}{\partial x}G\left(x,y,\lambda\right) -\frac{\partial^2}{\partial y\partial x}G\left(x,y,\lambda\right)\right]$, and thus our integrand will reduce to $\frac{\partial^2}{\partial^2 \lambda} \left[\frac{\partial}{\partial y}G\left(x,y,\lambda\right)\frac{\partial}{\partial x}G\left(x,y,\lambda\right)\right]e^{-G\left(x,y,\lambda\right)}$. Thus our integral is:
\begin{align*}
    &\int_{-\infty}^{\infty}\int_{-\infty}^{\infty}\frac{\partial^2}{\partial^2 \lambda} \left[\frac{\partial}{\partial y}G\left(x,y,\lambda\right)\frac{\partial}{\partial x}G\left(x,y,\lambda\right)\right]e^{-G\left(x,y,\lambda\right)}dxdy\\
    &= \int_{-\infty}^{\infty}\int_{-\infty}^{\infty}e^{-x}\phi\left(\lambda - \frac{\delta}{2\lambda}\right)\left[2e^{-x}\phi\left(\lambda - \frac{\delta}{2\lambda}\right)\left(1-\frac{\delta^2}{2\lambda^2}\right)\right.\\
    &\hspace{10mm}+e^{-x}\Phi\left(\lambda - \frac{\delta}{2\lambda}\right)\left\{-\lambda + \frac{\delta}{2\lambda}+\frac{\delta^2}{4\lambda^3}-\frac{\delta^3}{8\lambda^5}+\frac{\delta}{\lambda^3}\right\}\\
    &\hspace{10mm}+\left.e^{\delta-x}\Phi\left(\lambda + \frac{\delta}{2\lambda}\right)\left\{-\lambda - \frac{\delta}{2\lambda}+\frac{\delta^2}{4\lambda^3}+\frac{\delta^3}{8\lambda^5}-\frac{\delta}{\lambda^3}\right\}\right]e^{-e^{-x}\left(\Phi\left(\lambda - \frac{\delta}{2\lambda}\right)+e^{\delta}\Phi\left(\lambda + \frac{\delta}{2\lambda}\right)\right)}dxd\delta.
\end{align*}

Again, by factoring we see each term is bounded by Lemma~\ref{integral}. We move to the next term, $\frac{-\frac{\partial^4}{\partial y\partial x\partial^2\lambda}G\left(x,y,\lambda\right)}{\frac{\partial}{\partial y}G\left(x,y,\lambda\right)\frac{\partial}{\partial x}G\left(x,y,\lambda\right) -\frac{\partial^2}{\partial y\partial x}G\left(x,y,\lambda\right)}$. The form of the numerator is
\begin{align*}
    &\frac{\partial^4}{\partial y\partial x\partial^2\lambda}G\left(x,y,\lambda\right) = \frac{\partial}{\partial \lambda}\left[e^{-x}\phi\left(\lambda + \frac{y-x}{2\lambda}\right)\left(\frac{1}{2\lambda^2} + \frac{1}{2} - \frac{\left(y-x\right)^2}{8\lambda^4}\right)\right]\\
    & =e^{-x}\phi\left(\lambda + \frac{y-x}{2\lambda}\right)\left[\frac{-\lambda}{2}-\frac{1}{2\lambda}-\frac{1}{\lambda^3}+\frac{\left(y-x\right)^2}{4\lambda^3}+\frac{5\left(y-x\right)^2}{8\lambda^5}-\frac{\left(y-x\right)^4}{32\lambda^7}\right].
\end{align*}

Once again, the denominator of the term cancels with the same term in our density, and we see that each term is bounded by an integral of the form in Lemma~\ref{integral}. Finally, we move to the squared term: $\frac{\frac{\partial}{\partial \lambda} \left[\frac{\partial}{\partial y}G\left(x,y,\lambda\right)\frac{\partial}{\partial x}G\left(x,y,\lambda\right)\right] - \frac{\partial^3}{\partial y\partial x\partial\lambda}G\left(x,y,\lambda\right)^2}{\left(\frac{\partial}{\partial y}G\left(x,y,\lambda\right)\frac{\partial}{\partial x}G\left(x,y,\lambda\right) -\frac{\partial^2}{\partial y\partial x}G\left(x,y,\lambda\right)\right)^2}$. Again, one of the terms in the denominator cancels with that in the density, and we see:

\begin{align*}
    &0 \leq  \frac{\left(\frac{\partial}{\partial \lambda} \left[\frac{\partial}{\partial y}G\left(x,y,\lambda\right)\frac{\partial}{\partial x}G\left(x,y,\lambda\right)\right] - \frac{\partial^3}{\partial y\partial x\partial\lambda}G\left(x,y,\lambda\right)\right)^2}{\frac{\partial}{\partial y}G\left(x,y,\lambda\right)\frac{\partial}{\partial x}G\left(x,y,\lambda\right) -\frac{\partial^2}{\partial y\partial x}G\left(x,y,\lambda\right)} \\
    & = \frac{\left(e^{-x}\phi\left(\lambda + \frac{y-x}{2\lambda}\right)\left[e^{-x}\Phi\left(\lambda + \frac{y-x}{2\lambda}\right)\left(1-\frac{x-y}{2\lambda^2}\right) +e^{-y}\Phi\left(\lambda +\frac{x-y}{2\lambda}\right)\left(1-\frac{y-x}{2\lambda^2}\right)-\frac{1}{2\lambda^2}-\frac{1}{2}+\frac{\left(x-y\right)^2}{8\lambda^4}\right]\right)^2}{e^{-x}e^{-y}\Phi\left(\lambda + \frac{x-y}{2\lambda}\right)\Phi\left(\lambda + \frac{y-x}{2\lambda}\right) + e^{-x}\phi\left(\lambda + \frac{y-x}{2\lambda}\right)\frac{1}{2\lambda}}\\
    & \leq \frac{e^{-x}\phi\left(\lambda + \frac{y-x}{2\lambda}\right)^2\left(e^{-x}\Phi\left(\lambda + \frac{y-x}{2\lambda}\right)\left(1-\frac{x-y}{2\lambda^2}\right) +e^{-y}\Phi\left(\lambda +\frac{x-y}{2\lambda}\right)\left(1-\frac{y-x}{2\lambda^2}\right)-\frac{1}{2\lambda^2}-\frac{1}{2}+\frac{\left(x-y\right)^2}{8\lambda^4}\right)^2}{\phi\left(\lambda + \frac{y-x}{2\lambda}\right)\frac{1}{2\lambda}}\\
    & \leq 2\lambda e^{-x}\phi\left(\lambda + \frac{y-x}{2\lambda}\right)\left(e^{-x}\Phi\left(\lambda + \frac{y-x}{2\lambda}\right)\left(1-\frac{x-y}{2\lambda^2}\right) +e^{-y}\Phi\left(\lambda +\frac{x-y}{2\lambda}\right)\left(1-\frac{y-x}{2\lambda^2}\right)\right.\\
    &\hspace{10mm}-\left.\frac{1}{2\lambda^2}-\frac{1}{2}+\frac{\left(x-y\right)^2}{8\lambda^4}\right)^2.
\end{align*}

We see each of these terms is bounded by a form of the integrand in Lemma~\ref{integral}. This was the final term of the integrand, and thus the Fisher information is finite.

\qed

We note that the integrand in our definition of the Fisher Information is continuous in $\lambda$ on the interior of the parameter space. Thus by the Leibniz Integral Rule, we may pass a derivative under the integral and, by the usual argument, we may represent the Fisher Information as the negation of the second moment of the first derivative of the log-likelihood function with respect to $\lambda$. This is to say that the value of the Fisher Information is indeed positive.

\section{Estimation Under the Misspecified Model}
\textcolor{black}{We now move from maximum likelihood estimation under the correctly specified model to that under the model misspecified by BM. Namely, let $(X_i,Y_i)$ for $i \in \{i,...,n\}$ be a sample of i.i.d. standard bivariate normal random variables, with correlation parameter $\rho_m$. Then, we are interested in the properties of the MLE
\begin{equation} \label{MLEdef}
    \hat{\lambda}_n := \argmax_{\lambda \in (0,\infty)}\prod_{j=1}^{k_n} h_{\lambda}(M_{X,j}, M_{Y,j}),
\end{equation}
where for $j \in \{i,...,k_n\}$ and $Z \in {X,Y}$,
\begin{equation*}
    M_{Z,j} = \max\{b_mZ_{j} - b_m^2;j=jm_n,...,(j+1)m_n-1\},
\end{equation*} are finite, where $b_m$ is defined by $b_m = m\phi(b_m)$, $\phi$ being the standard normal density, and $h_\lambda$ being the H\"usler-Reiss density with parameter $\lambda$. We note that, at the moment, we are not assured that such a value necessarily exists. The goal of this section is to derive asymptotic properties of this estimator under the block maxima method. The final result takes the following form:}

\begin{theorem}\label{MLE Form}
    \textcolor{black}{Define the MLE by (\ref{MLEdef}). Assume that the true value of $\lambda \in (0,\infty)$. Assume that $k= k_n \rightarrow \infty$ and $m = m_n \rightarrow \infty$, and are defined such that 
    \begin{equation}\label{limitAssumptions}
    L_1:=\lim_{m\rightarrow\infty}\frac{\sqrt{k}}{b_m^2}, \hspace{5mm} L_2:=\lim_{m\rightarrow\infty}\sqrt{k}\left(\lambda - b_m\sqrt{\frac{1-\rho_m}{1+\rho_m}}\right)
\end{equation}exist and are finite. There exists a sequence of estimators $\hat{\lambda}_n$ such that 
    \begin{equation*}
        \lim_{n\rightarrow\infty}\mathbb{P}\left(\hat{\lambda}_n \text{ is an MLE}\right) = 1
    \end{equation*} and
    \begin{equation*}
        \sqrt{k}\left(\hat{\lambda}_n-\lambda\right) \xrightarrow{d} N\left(I_\lambda^{-1}A,I_\lambda^{-1}\right)
    \end{equation*}
    where $A \in \mathbb{R}$ is the limiting bias defined in Lemma~\ref{bias-variance theorem}.}
\end{theorem}
\textcolor{black}{
\begin{remark}
    Note that $\lim_{m\rightarrow\infty}\frac{\sqrt{k}}{b_m^2} = \lim_{m\rightarrow\infty}\frac{\sqrt{n/m}}{2\ln(m)}$, as $n \approx k m$. As a result, when $L_1>0$ we have $n \approx L_1^2m\ln(m)^2$ and thus $\frac{\ln(n)}{\ln(m)} \to 1$. However, we allow for the case that $L_1 = 0$, and thus we only have $n = O(m\ln(m)^2)$ in general, which gives us a first condition on the balance of $m$ and $k$.
\end{remark}}\textcolor{black}{It is from this result that we may derive results on the best practices for choosing $k$ and $m$ and give an explicit form for the asymptotic bias of the MLE. We note that our notation differs from some others such as \textcite{huslerreiss89}, in order to preserve the convention that $n$ represents overall sample size. As such, $\rho_m$ and thus the value of $\lambda$ depends on the block size $m$ rather than $n$. This section includes a few lemmas necessary to prove Theorem \ref{MLE Form} and a short proof of theorem itself. The form of the asymptotic bias term $A$ together with some auxiliary results are relegated to the appendix. We begin with a proof of the convergence of the Fisher Informations.}
\begin{redlemma}\label{Convergence of expectations}
Let ($X_i$, $Y_i$) for $i \in \left\{1,..,m\right\}$ be independent and identically distributed according to a standard bivariate normal distribution with correlation parameter $\rho_m$ such that $\ln\left(m\right)\left(1-\rho_m\right) \rightarrow \lambda^2 \in \left(0, \infty\right)$. Define $l''\left(x,y,\lambda\right)$ to be the second derivative with respect to $\lambda$ of the log-likelihood function of a H\"usler-Reiss distribution. Namely, $l''\left(x,y,\lambda\right) := \frac{\partial^2}{\partial^2\lambda}\left[\ln\left(h_{\lambda}\left(x,y\right)\right)\right]$, where $h_{\lambda}\left(x,y\right)$ is the density function of the H\"usler-Reiss distribution with parameter $\lambda$. Then the expectation of $l''\left(x,y,\lambda\right)$ with respect to the density of the maxima of the properly scaled $\left(X_i,Y_i\right)$ converges to $I_\lambda$, the finite Fisher Information of the H\"usler-Reiss distribution, which depends only on the value of $\lambda$.
\end{redlemma}

\noindent \textit{Proof.}
Following \textcites{huslerreiss89} in their choice of scaling constants, we choose $b_m$, defined by $b_m = m\phi(b_m)$ for both the location and scale constants. The first step of the proof is to show that the density of $\left(b_m\max\left(X_i\right)-b_m^2,b_m\max\left(Y_i\right)-b_m^2\right)$ converges to $h_{\lambda}\left(x,y\right)$. Define $u_m\left(x\right) = \frac{x}{b_m} + b_m$. Recall that $\lambda \in \left(0,\infty\right)$ implies that $\rho_m \rightarrow 1$. \textcolor{black}{It will be useful moving forward to have defined the equation
\begin{equation}
    q_m(x,y) = \frac{u_m(x)-\rho_m u_m(y)}{\sqrt{1-\rho_m^2}}.
\end{equation}}Recall from Equations 2.11 and 2.12 of \textcites{huslerreiss89} that, without loss of generality, $q_m(x,y) \rightarrow \lambda + \frac{x-y}{2\lambda}$. Define $\Phi_{\rho_m}\left(x,y\right)$ and $\phi_{\rho_m}\left(x,y\right)$ as the CDF (cumulative distribution function) and PDF (probability density function) of a standard normal bivariate distribution with correlation parameter $\rho_m$. Throughout the proof we will utilize the following useful equality: $\phi_{\rho_m}\left(x,y\right) = \phi\left(\frac{x-\rho_my}{\sqrt{1-\rho_m^2}}\right)\phi\left(y\right)\left(1-\rho_m^2\right)^{-1/2}$. We begin by differentiating the CDF:
\begin{align*}
    &\frac{\partial^2}{\partial x\partial y}\Phi_{\rho_m}\left(u_m\left(x\right),u_m\left(y\right)\right)^m \\
    &= \Phi_{\rho_m}\left(u_m\left(x\right),u_m\left(y\right)\right)^{m-2}\left[\frac{m-1}{m}e^{-y}e^{-\frac{y^2}{2b_m^2}}\Phi\left(\textcolor{black}{q_m(x,y)}\right)e^{-x}e^{-\frac{x^2}{2b_m^2}}\Phi\left(\textcolor{black}{q_m(y,x)}\right)\right.\\
    &\left.\hspace{10mm}+ \Phi_{\rho_m}\left(u_m\left(x\right),u_m\left(y\right)\right)e^{-x}e^{-\frac{x^2}{2b_m^2}}\phi\left(\textcolor{black}{q_m(y,x)}\right)\frac{1}{b_m\sqrt{1-\rho_m^2}}\right].
\end{align*}
We note that $\Phi_{\rho_m}\left(u_m\left(x\right),u_m\left(y\right)\right),e^{-\frac{x^2}{2b_m^2}},e^{-\frac{y^2}{2b_m^2}},\frac{m-1}{m} \rightarrow 1$, $\Phi_{\rho_m}\left(u_m\left(x\right),u_m\left(y\right)\right)^{m-2} \rightarrow H_{\lambda}\left(x,y\right)$ and $b_m\sqrt{1-\rho_m^2} \rightarrow 2\lambda$ as $m \rightarrow \infty$. This gives us:
\begin{align*}
    &\frac{\partial^2}{\partial x\partial y}\Phi_{\rho_m}\left(u_m\left(x\right),u_m\left(y\right)\right)^m \\
    &= \Phi_{\rho_m}\left(u_m\left(x\right),u_m\left(y\right)\right)^{m-2}\left[\frac{m-1}{m}e^{-y}e^{-\frac{y^2}{2b_m^2}}\Phi\left(\textcolor{black}{q_m(x,y)}\right)e^{-x}e^{-\frac{x^2}{2b_m^2}}\Phi\left(\textcolor{black}{q_m(y,x)}\right)\right.\\
    &\hspace{10mm}\left.+ \Phi_{\rho_m}\left(u_m\left(x\right),u_m\left(y\right)\right)e^{-x}e^{-\frac{x^2}{2b_m^2}}\phi\left(\textcolor{black}{q_m(y,x)}\right)\frac{1}{b_m\sqrt{1-\rho_m^2}}\right] \\
    &\rightarrow H_{\lambda}\left(x,y\right)\left[e^{-x}e^{-y}\Phi\left(\lambda + \frac{x-y}{2\lambda}\right)\Phi\left(\lambda + \frac{y-x}{2\lambda}\right) + e^{-x}\phi\left(\lambda + \frac{y-x}{2\lambda}\right)\frac{1}{2\lambda}\right] = h_\lambda(x,y)
\end{align*}
This limiting density is exactly the density of a H\"usler-Reiss distribution with parameter $\lambda$ derived in Theorem~\ref{MLE existence}, therefore the densities converge pointwise. Thus, the limiting expectation is exactly the Fisher Information which was proven to be finite in Theorem~\ref{MLE existence}. To show that these expectations converge we will use the Dominated Convergence Theorem (DCT). We do this in two steps: first, we show that 
\begin{equation}\label{goal1}
    \iint_{\mathbb{R}^2} l''\left(x,y,\lambda\right)\frac{\partial^2}{\partial x\partial y}\Phi_{\rho_m}\left(u_m\left(x\right),u_m\left(y\right)\right)^m1_{x \geq -b_m^2}dxdy \rightarrow I_\lambda,
\end{equation}
then, we show that 
\begin{equation}\label{goal2}
    \iint_{\mathbb{R}^2} l''\left(x,y,\lambda\right)\frac{\partial^2}{\partial x\partial y}\Phi_{\rho_m}\left(u_m\left(x\right),u_m\left(y\right)\right)^m1_{x \leq -b_m^2} \rightarrow 0,
\end{equation}
which in turn proves the theorem. To show (\ref{goal1}), we want to show there exists a function $f\left(x,y\right)$ that is independent of $m$ such that $|l''\left(x,y,\lambda\right)\frac{\partial^2}{\partial x\partial y}\Phi_{\rho_m}\left(u_m\left(x\right),u_m\left(y\right)\right)^m1_{x \geq -b_m^2}| \leq f\left(x,y\right)$ and $\iint_{\mathbb{R}^2}f\left(x,y\right)1_{x \geq -b_m^2}dxdy < \infty$. We assume

\begin{equation}\label{form 1}
    \int_{\mathbb{R}}|l''\left(x,y,\lambda\right)\frac{\partial^2}{\partial x\partial y}\Phi_{\rho_m}\left(u_m\left(x\right),u_m\left(y\right)\right)^m|dy \leq e^{-kx}\Phi\left(u_m\left(x\right)\right)^{m-2}e^{-\frac{x^2}{2b_m^2}}
\end{equation}
for some $k > 0$. With this assumption, we show (\ref{goal1}) and (\ref{goal2}), and then go back to prove the assumption. Assuming the form $e^{-kx}\Phi\left(u_m\left(x\right)\right)^{m-2}e^{-\frac{x^2}{2b_m^2}}$, we see
\begin{align*}
    &\int_0^{\infty}e^{-kx}\Phi\left(u_m\left(x\right)\right)^{m-2}e^{-\frac{x^2}{2b_m^2}}dx \leq \int_{0}^{\infty}e^{-kx}dx = \frac{1}{k}.
\end{align*}
Consider the integral over $\left(-kb_m,0\right)$ for some fixed $k > 1$. We see
\begin{align*}
    \int_{-kb_m}^{0}e^{-kx}\Phi\left(u_m\left(x\right)\right)^{m-2}e^{-\frac{x^2}{2b_m^2}}dx
    \leq & \int_{-kb_m}^{0}e^{-kx}\Phi\left(u_m\left(x\right)\right)^{m-2}dx\\
    \leq & \int_{-kb_m}^{0}e^{-kx}\exp\left[-\frac{m-2}{m}\int_{x}^{\infty}e^{-t}e^{-\frac{1}{2}\left(\frac{t}{b_m}\right)^2}dt\right]dx\\
    \leq &\int_{-kb_m}^{0}e^{-kx}\exp\left[-\frac{1}{3}e^{-\frac{k^2}{2}}\int_{x}^{0}e^{-t}dt\right]dx\\
    \leq &e^{\frac{1}{3}e^{-\frac{k^2}{2}}}\int_{-kb_m}^{0}e^{-kx}e^{-\frac{1}{3}e^{-\frac{k^2}{2}}e^{-x}}dx.
\end{align*}
We see that the remaining integral is of the integrable form in Lemma~\ref{integral} with $f\left(x-y\right) = \frac{1}{3}e^{-\frac{k^2}{2}}$. Thus we see that the integral is bounded by a finite value over $\left[-kb_m,\infty\right)$. Now we need to show that the integral over $\left[-\infty,-kb_m\right]$ goes to zero as $m \rightarrow \infty$. Consider the integral over $\left(-\infty,-b_m^2\right)$:
\begin{align*}
    \int_{-\infty}^{-b_m^2}e^{-kx}\Phi\left(u_m\left(x\right)\right)^{m-2}e^{-\frac{x^2}{2b_m^2}}dx
    = & \sqrt{2\pi}e^{\frac{1}{2}k^2b_m^2}\int_{-\infty}^{-b_m^2}\phi\left(\frac{x}{b_m}+kb_m\right)\Phi\left(u_m\left(x\right)\right)^{m-2}dx\\
    \leq & \sqrt{2\pi}e^{\frac{1}{2}k^2b_m^2}\left(\frac{1}{2}\right)^{m-2}\int_{-\infty}^{-b_m^2}\phi\left(\frac{x}{b_m}+kb_m\right)dx\\
    \leq & \sqrt{2\pi}b_m m^{k^2}\left(\frac{1}{2}\right)^{m-2} \rightarrow 0.
\end{align*}
Finally, we look at the integral over the region $\left(-b_m^2, -kb_m\right)$: 
\begin{align*}
    &\int_{-b_m^2}^{-kb_m}e^{-kx}e^{-\frac{1}{2}\frac{x^2}{b_m^2}}\Phi\left(u_m\left(x\right)\right)^{m-2}dx.
\end{align*}
Taking the derivative of the integrand we see:
\begin{equation*}
    \frac{\partial}{\partial x}e^{-kx}e^{-\frac{1}{2}\frac{x^2}{b_m^2}}\Phi\left(u_m\left(x\right)\right)^{m-2} \textcolor{black}{\geq} e^{-kx}e^{-\frac{1}{2}\frac{x^2}{b_m^2}}\Phi\left(u_m\left(x\right)\right)^{m-3}\left[-\Phi\left(u_m\left(x\right)\right)\left(k + \frac{x}{b_m}\right) + \phi\left(\frac{x}{b_m}+b_m\right)/b_m\right].
\end{equation*}
This is surely positive for $x \leq -kb_m$. Let $m$ be sufficiently large such that $\frac{b_m-2}{b_m-1} > \frac{1}{2}$ Then, returning to our integral, we see:
\begin{align*}
    \int_{-b_m^2}^{-kb_m}e^{-kx}e^{-\frac{1}{2}\frac{x^2}{b_m^2}}\Phi\left(u_m\left(x\right)\right)^{m-2}dx
    &\leq b_m^2e^{k^2b_m}e^{-\frac{k^2}{2}}\Phi\left(b_m-k\right)^{m-2}\\
    &\leq b_m^2e^{k^2b_m}e^{-\frac{k^2}{2}}\Phi\left(b_m-1\right)^{m-2}\\
    &\leq e^{-\frac{k^2}{2}}e^{\left(k^2+1\right)b_m}\left(1-\left(1-\Phi\left(b_m-1\right)\right)\right)^{m-2}\\
\end{align*}\textcolor{black}{We continue by using the inequality $1 - \Phi(x) > \frac{\phi(x)}{x}(1-\frac{1}{x^2}$ for $x>0$. This gives us} 
\begin{align*}
    \textcolor{black}{\int_{-b_m^2}^{-kb_m}e^{-kx}e^{-\frac{1}{2}\frac{x^2}{b_m^2}}\Phi\left(u_m\left(x\right)\right)^{m-2}dx}&\leq e^{-\frac{k^2}{2}}e^{\left(k^2+1\right)b_m}\left(1-\frac{\phi\left(b_m-1\right)}{b_m-1}\left(1-\frac{1}{\left(b_m-1\right)^2}\right)\right)^{m-2}\\
    &= e^{-\frac{k^2}{2}}e^{\left(k^2+1\right)b_m}\left(1-\frac{1}{2m}e^{b_m-\frac{1}{2}}\right)^{-2}e^{-m\sum_{i=1}^\infty\frac{\left(\frac{1}{2}e^{b_m-\frac{1}{2}}\right)^i}{im^i}}\\
    &\leq e^{-\frac{k^2}{2}}e^{\left(k^2+1\right)b_m}\left(1-\frac{1}{2m}e^{b_m-\frac{1}{2}}\right)^{-2}e^{-\frac{1}{2}e^{b_m-\frac{1}{2}}}\rightarrow0.
\end{align*}

All that is left to prove is that (\ref{form 1}) is satisfied. To show this, we break $l''\left(x,y,\lambda\right)$ into its constituent terms and bound their absolute values. Recall the form of $l''\left(x,y,\lambda\right)$ as given in (\ref{second derivative form}) and made explicit in the following calculations. We see that the absolute value of each of these terms may be bounded above by a term of the form:
\begin{equation}\label{boundform}
    Ce^{-k_ix}e^{-j_iy}P_i\left(x-y\right),
\end{equation}
where $k_i, j_i$ are non-negative integers and $P_i\left(t\right)$ is a polynomial in $t$ for each $i \in \left\{1,2,3,4\right\}$. This variable $i$ is just an indexing variable which will not affect the arguments below. For brevity, we drop the subscript in the following analysis. 
Consider the integral of one of these terms.
\begin{align}\label{split}
    &\iint_{\mathbb{R}^2}Ce^{-kx}e^{-jy}P\left(x-y\right)\frac{\partial^2}{\partial x\partial y}\Phi_{\rho_m}\left(u_m\left(x\right),u_m\left(y\right)\right)^mdxdy\\
    &=C\iint_{\mathbb{R}^2}e^{-kx}e^{-jy}P\left(x-y\right)\Phi_{\rho_m}\left(u_m\left(x\right),u_m\left(y\right)\right)^{m-2}\nonumber\\
    &\hspace{10mm}\left[\frac{m-1}{m}e^{-y}e^{-\frac{y^2}{2b_m^2}}\Phi\left(\textcolor{black}{q_m(x,y)}\right)\nonumber\right.e^{-x}e^{-\frac{x^2}{2b_m^2}}\Phi\left(\textcolor{black}{q_m(y,x)}\right)\\
    &\hspace{10mm} +\left.\Phi_{\rho_m}\left(u_m\left(x\right),u_m\left(y\right)\right)e^{-x}e^{-\frac{x^2}{2b_m^2}}\phi\left(\textcolor{black}{q_m(y,x)}\right)\frac{1}{b_m\sqrt{1-\rho_m^2}}\right]dxdy\nonumber\\
     \leq &C\iint_{\mathbb{R}^2}e^{-\left(k+1\right)x}e^{-jy}P\left(x-y\right)\Phi\left(u_m\left(x\right)\right)^{m-2}e^{-\frac{x^2}{2b_m^2}}\nonumber\\
     &\hspace{10mm}\left[e^{-y}e^{-\frac{y^2}{2b_m^2}}\Phi\left(\textcolor{black}{q_m(x,y)}\right)\Phi\left(\textcolor{black}{q_m(y,x)}\right)+\phi\left(\textcolor{black}{q_m(y,x)}\right)\frac{1}{b_m\sqrt{1-\rho_m^2}}\right]dxdy.\nonumber
\end{align}

For further analysis, we break the integral across its sum. We also assume that $m$ is sufficiently large such that $\rho_m > 0$. To begin the analysis, we first notice that $b_m\sqrt{1-\rho_m^2}$ \textcolor{black}{$\rightarrow 2\lambda \neq 0$ is a convergent sequence whose inverse will also be convergent and thus bounded. Define $M_1$ to be one such upper bound on the inverse of the sequence. Likewise, let $M_2$ be an upper bound of $b_m\sqrt{1-\rho_m^2}$.} There must exist $l_1,l_2,C$ such that $P\left(\delta\right) \leq C\left(e^{-l_1\delta}+e^{l_2\delta}\right)$. Using this bound and the transformation of variables $\delta = x-y$, we analyze the second term:
\begin{align*}
    &\iint_{\mathbb{R}^2}e^{-\left(k+1\right)x}e^{-jy}P\left(x-y\right)\Phi\left(u_m\left(x\right)\right)^{m-2}e^{-\frac{x^2}{2b_m^2}}\left(\phi\left(\textcolor{black}{q_m(y,x)}\right)\frac{1}{b_m\sqrt{1-\rho_m^2}}\right)dxdy\\
    & \leq M_1\int e^{-\left(k+j+1\right)x}\Phi\left(u_m\left(x\right)\right)^{m-2}e^{-\frac{x^2}{2b_m^2}}\left(\int e^{j\delta}P\left(\delta\right)\phi\left(\frac{\delta-\left(1-\rho_m\right)\left(b_m^2-x\right)}{b_m\sqrt{1-\rho_m^2}}\right)d\delta\right)dx\\
    & \leq CM_1\int e^{-\left(k+j+1\right)x}\Phi\left(u_m\left(x\right)\right)^{m-2}e^{-\frac{x^2}{2b_m^2}} \left(\int \left(e^{j-l_1\delta}+e^{\left(l_2+j\right)\delta}\right)\phi\left(\frac{\delta-\left(1-\rho_m\right)\left(b_m^2-x\right)}{b_m^2\sqrt{1-\rho_m^2}}\right)d\delta\right)dx\\
    & = CM_1\int e^{-\left(k+j+1\right)x}\Phi\left(u_m\left(x\right)\right)^{m-2}e^{-\frac{x^2}{2b_m^2}}\\
    &\hspace{10mm}\left(e^{\left(j-l_1\right)\left(1-\rho_m\right)\left(b_m^2-x\right) + \left(j-l_1\right)^2/2b_m\left(1-\rho_m^2\right)} + e^{\left(l_2+j\right)\left(1-\rho_m\right)\left(b_m^2-x\right) + \left(l_2+j\right)^2/2b_m\left(1-\rho_m^2\right)}\right)dx\\
    & \leq CM_1e^{\max\left(|j-l_1|,
    l_2+j\right)M_1^2+\max\left(|j-l_1|,l_2+j\right)^2M_2}\\
    &\hspace{10mm}\int e^{-\left(k+j+1\right)x}\Phi\left(u_m\left(x\right)\right)^{m-2}e^{-\frac{x^2}{2b_m^2}}\left(e^{-\left(j-l_1\right)\left(1-\rho_m\right)x} + e^{-\left(l_2+j\right)\left(1-\rho_m\right)x}\right)dx,
\end{align*} 
where the transition from the third to fourth line is dependent on the moment generating function of the normal distribution. But for sufficiently large $m$, we know that $-\left(l_2+j\right)\left(1-\rho_m\right) - \left(k+j+1\right) < 0$ and $-\left(j-l_1\right)\left(1-\rho_m\right) - \left(k+j+1\right) < 0$ as $\rho_m \rightarrow 1$. Thus, for sufficiently large $m$ we are left with the integral:
\begin{equation*}
    \int \left(e^{-s_1x} + e^{-s_2x}\right)\Phi\left(u_m\left(x\right)\right)^{m-2}e^{-\frac{x^2}{2b_m^2}}dx.
\end{equation*}
This works for any term in our integrand. Throughout this process, all inequalities have come from bounding the integrand itself, not the integral. Thus, noting that this last integral is of the form in (\ref{form 1}), we see that these terms satisfy all requirements for this theorem. We still must prove the same requirements are met for the first term in the sum in (\ref{split}). 

It is easy to check $\Phi\left(x\right) \leq 2\phi\left(x\right)$ for $x<0$. This leads us to $\Phi(\textcolor{black}{q_m(x,x-\delta)}) \leq 2\phi(\textcolor{black}{q_m(x,x-\delta)})$ when $\delta \textcolor{black}{<} -\frac{1-\rho_m}{\rho_m}\left(x+b_m^2\right)$ and $\Phi(\textcolor{black}{q_m(x-\delta,x)}) \leq 2\phi(\textcolor{black}{q_m(x-\delta,x)})$ when $\delta \textcolor{black}{>} \left(1-\rho_m\right)\left(x+b_m^2\right)$. For the purpose of the following analysis, we assume $m$ is sufficiently large such that $\rho_m > \frac{1}{2}$. Recalling that $\left(1-\rho_m\right)b_m^2 \rightarrow 2\lambda^2$, we note that this \textcolor{black}{sequence must be bounded above, and denote such upper bound by $M_3$}. Then, we can relax our conditions to $\delta > \textcolor{black}{\left(1-\rho_m\right)}x + M_3$ and $\delta < -\textcolor{black}{\frac{1-\rho_m}{\rho_m}}x-2M_3$, \textcolor{black}{respectively}. With this, we may finally make the calculations:
\begin{align*}
    &\iint_{\mathbb{R}^2}Ce^{-kx}e^{-jy}P\left(x-y\right)\Phi_{\rho_m}\left(u_m\left(x\right),u_m\left(y\right)\right)^{m-2}e^{-y-x}e^{-\frac{y^2+x^2}{2b_m^2}}\Phi\left(\textcolor{black}{q_m(x,y)}\right)\Phi\left(\textcolor{black}{q_m(y,x)}\right)dxdy\\
    \leq &\int\Phi\left(u_m\left(x\right)\right)^{m-2}e^{-\left(2+k+j\right)x}e^{-\frac{x^2}{2b_m^2}}\left(\int e^{(j+1)\delta}P\left(\delta\right)\Phi\left(\textcolor{black}{q_m(x,x-\delta)}\right)\Phi\left(\textcolor{black}{q_m(x-\delta,x)}\right)d\delta\right)dx\\
    \leq &\int\Phi\left(u_m\left(x\right)\right)^{m-2}e^{-\left(2+k+j\right)x}e^{-\frac{x^2}{2b_m^2}}\left(\int_{\left(1-\rho_m\right)x+M_{\textcolor{black}{3}}}^{\infty}2e^{\left(j+1\right)\delta}P\left(\delta\right)\phi\left(\textcolor{black}{q_m(x,x-\delta)}\right)d\delta\right.\\
    +& \left.\int^{\left(1-\rho_m\right)x+M_{\textcolor{black}{3}}}_{-\frac{1-\rho_m}{\rho_m}x-2\textcolor{black}{M_3}}e^{\left(j+1\right)\delta}P\left(\delta\right)d\delta + \int_{-\infty}^{-\frac{1-\rho_m}{\rho_m}x-2M_{\textcolor{black}{3}}}2e^{\left(j+1\right)\delta}P\left(\delta\right)\phi\left(\textcolor{black}{q_m(x-\delta,x)}\right)d\delta\right)dx.
\end{align*}
Note that the terms including $\int_{-\infty}^{-\frac{1-\rho_m}{\rho_m}x-2M_{\textcolor{black}{3}}}2e^{\left(j+1\right)\delta}P\left(\delta\right)\phi\left(\textcolor{black}{q_m(x-\delta,x)}\right)d\delta$ and \\$\int_{\left(1-\rho_m\right)x+M_{\textcolor{black}{3}}}^{\infty}2e^{\left(j+1\right)\delta}P\left(\delta\right)\phi\left(\textcolor{black}{q_m(x,x-\delta)}\right)d\delta$ are of the same form shown to satisfy our conditions in our analysis of the first term in (\ref{split}). Seeing out the calculation of the remaining term and utilizing $P\left(\delta\right) \leq C_1e^{k_1\delta} + C_2e^{-k_2\delta}$:

\begin{align*}
    &\int^{\left(1-\rho_m\right)x+M_{\textcolor{black}{3}}}_{-\frac{1-\rho_m}{\rho_m}x-2M_{\textcolor{black}{3}}}e^{\left(j+1\right)\delta}P\left(\delta\right)d\delta \\
    \leq &\int^{\left(1-\rho_m\right)x+M_{\textcolor{black}{3}}}_{-\frac{1-\rho_m}{\rho_m}x-2M_{\textcolor{black}{3}}}C_1e^{\left(k_1+j+1\right)\delta} + C_2e^{\left(j+1-k_2\delta\right)}d\delta\\
    = & \frac{C_1}{k_1+j+1}\left[e^{\left(k_1+j+1\right)\left(1-\rho_m\right)x + M_{\textcolor{black}{3}}\left(k_1+j+1\right)}-e^{-\left(k_1+j_1+1\right)\frac{\left(1-\rho_m\right)}{\rho_m}x - 2M_{\textcolor{black}{3}}\left(k_1+j+1\right)}\right]\\
   &\hspace{10mm}+ \frac{C_2}{j+1-k_2}\left[e^{\left(j+1-k_2\right)\left(1-\rho_m\right)x + M_{\textcolor{black}{3}}\left(j+1-k_2\right)}-e^{\left(k_2+j+1\right)\frac{\left(1-\rho_m\right)}{\rho_m}x - 2M_{\textcolor{black}{3}}\left(j+1-k_2\right)}\right].
\end{align*}

Once again, we see that this leaves the integral in the form of (\ref{form 1}). Thus all three of the final terms can be bounded and by Dominated Convergence Theorem, the expectations converge. \qed

\sloppy We note that this immediately implies that the variance of this term is convergent and thus also bounded. Since each term in $\frac{\partial^2}{\partial^2 \lambda}\ln\left(h_{\lambda}\left(x,y\right)\right)$ was treated as $Ce^{-k_ix-j_iy}P_i\left(x-y\right)$, taking the square of or product of any two terms yields the same form. Thus, the same argument above applies to the variance. 

\textcolor{black}{With this result, we can begin to build up the structure necessary to prove Theorem~\ref{MLE Form}. Our approach will to be to define a localized log-likelihood function $\tilde{L}$ and analyze its Taylor series expansion as the sample size $n$ grows. We begin with the following lemma:}

\begin{redlemma}\label{op theorem}
\textcolor{black}{Let k, m satisfy the same conditions as in Theorem \ref{MLE Form}. Let $x_i,y_i$ for $i =1,..,k$ to be the scaled maxima of $m$ bivariate normal random variables with mean 0, unit variance, and correlation $\rho_m$ such that the limiting H\"usler-Reiss dependence parameter $\lambda \in (0,
\infty)$.} Denote $\tilde{L}\left(h\right) = \displaystyle\sum_{i=1}^k\ln\left(h_{\lambda+\frac{h}{\sqrt{k}}}\left(x_i,y_i\right)\right)$ with $h \in H_k$ where $H_k$ is a ball of radius $r_k = o\left(\sqrt{k}\right)$, \textcolor{black}{centered at the origin}. Then, we have:
\begin{equation}\label{Op}
    \frac{\partial^2}{\partial^2 h}\tilde{L}\left(h\right) - I_\lambda = o_p\left(1\right).
\end{equation}
\end{redlemma}

\noindent \textit{Proof.}
Denote $\frac{\partial^2}{\partial^2 \lambda}\ln\left(h_{\lambda+\frac{h}{\sqrt{k}}}\left(x_i,y_i\right)\right) = l''\left(h,x,y\right)$ and note that since $h$ must grow slower than $\sqrt{k}$, $\lambda + \frac{h}{\sqrt{k}}$ stays within an open interval of $\lambda$. \textcolor{black}{Although not strictly necessary, the condition that $H_k$ be centered at the origin ensures that $h = 0$ contained in our domain.} Throughout this proof we will need to take expectations with respect to the limiting H\"usler-Reiss distribution and with respect to the distribution of the maxima of scaled bivariate normal random variables. We will denote the former $\mathbb{E}_\lambda$ and the latter $\mathbb{E}_{\phi_m}$. Applying the triangle inequality, we break the difference into three other differences:
\begin{align*}
    \left|\frac{\partial^2}{\partial^2 h}\tilde{L}\left(h\right) - I_\lambda\right| &=  \left|\frac{1}{k}\sum_{i=1}^k\frac{\partial^2}{\partial^2 \lambda}\ln\left(\frac{\partial^2}{\partial x\partial y}H\left(x_i,y_i\right)_{\lambda+\frac{h}{\sqrt{k}}}\right) - \mathbb{E}_{\phi_m}\left(l''\left(h,X,Y\right)\right)\right. \\
    & \hspace{10mm}\left.+ \mathbb{E}_{\phi_m}\left(l''\left(h,X,Y\right)\right)-\mathbb{E}_{\phi_m}\left(l''\left(0,X,Y\right)\right) + \mathbb{E}_{\phi_m}\left(l''\left(0,X,Y\right)\right) - I_\lambda\right|\\
    & \leq  \left|\frac{1}{k}\sum_{i=1}^k\frac{\partial^2}{\partial^2 \lambda}\ln\left(\frac{\partial^2}{\partial x\partial y}H\left(x_i,y_i\right)_{\lambda+\frac{h}{\sqrt{k}}}\right) - \mathbb{E}_{\phi_m}\left(l''\left(h,X,Y\right)\right)\right| \\
    & \hspace{10mm}+  \left|\mathbb{E}_{\phi_m}\left(l''\left(h,X,Y\right)\right) -\mathbb{E}_{\phi_m}\left(l''\left(0,X,Y\right)\right)\right|\\
    & \hspace{10mm}+  \left|\mathbb{E}_{\phi_m}\left(l''\left(0,X,Y\right)\right) - I_\lambda\right|.
\end{align*}
Denote $|\mathbb{E}_{\phi_m}\left(l''\left(0,X,Y\right)\right) - I_\lambda| = a_m$. As shown in \textcolor{black}{Lemma}~\ref{Convergence of expectations}, $\mathbb{E}_{\phi_m}\left(l''\left(0,X,Y\right)\right) \rightarrow I_\lambda$, and thus $a_m \rightarrow 0$. Denote $|\mathbb{E}_{\phi_m}\left(l''\left(h,X,Y\right)\right) -\mathbb{E}_{\phi_m}\left(l''\left(0,X,Y\right)\right)| = c_m$. Since $h \in H_k$, we know that $\lambda + \frac{h}{\sqrt{k}} \rightarrow \lambda$. Consider, then,

\begin{align*}
    &\mathbb{E}_{\phi_m}\left(l''\left(h,X,Y\right)\right) -\mathbb{E}_{\phi_m}\left(l''\left(0,X,Y\right)\right)\\
    &= \iint_{\mathbb{R}^2}\left(l''\left(h,x,y\right) - l''\left(0,x,y\right)\right)\frac{\partial^2}{\partial x\partial y}\Phi_{\rho_m}\left(u_m\left(x\right),u_m\left(y\right)\right)^m dxdy.
\end{align*}
Since $l$ is continuous, this integrand converges to 0 as $n \rightarrow \infty$. We also note that $c_m \leq \mathbb{E}_{\phi_m}\left|l''\left(h,X,Y\right)\right| + \mathbb{E}_{\phi_m}\left|l''\left(0,X,Y\right)\right|$.
The second of these terms can be bounded by DCT as per \textcolor{black}{Lemma}~\ref{Convergence of expectations}, and \textcolor{black}{$c_m \to 0$}.  To take care of the first term, we need to show that we can once again use DCT to show the convergence of this integral, even when $\lambda$ is not exactly equal to its true value, but instead converges to it. \textcolor{black}{This reduces to showing that we may} bound each term by the form of (\ref{boundform}). We notice that any $\Phi$ or $\phi$ is immediately negated as these are bounded functions, and do not appear in the form of the bound. Thus, the only time $\lambda + \frac{h}{\sqrt{k}}$ appears is in the polynomial $P_i(x-y)$. Since $\lambda + \frac{h}{\sqrt{k}}$ converges to a non-zero quantity, for sufficiently large $k$ both $\lambda + \frac{h}{\sqrt{k}}$ and $(\lambda + \frac{h}{\sqrt{k}})^{-1}$ will have finite suprema. Thus, wherever $\lambda + \frac{h}{\sqrt{k}}$ appears in this bound, the we may bound it by a value which does not depend on $m$ and use the same form in (\ref{boundform}). Thus, the same argument present in \textcolor{black}{Lemma}~\ref{Convergence of expectations} applies, and by DCT the quantity will converge. Fix $\epsilon > 0$. Using Chebyshev's Inequality, we can now show:

\begin{align*}
    &\mathbb{P}\left(\left|\frac{\partial^2}{\partial^2 h}\tilde{L}\left(h\right) - I_\lambda\right| > \epsilon\right) \\
    &\leq \mathbb{P}\left(\left|\frac{1}{k}\sum_{i=1}^k\frac{\partial^2}{\partial^2 \lambda}\ln\left(\frac{\partial^2}{\partial x\partial y}H\left(X_i,Y_i\right)_{\lambda+\frac{h}{\sqrt{k}}}\right) - \mathbb{E}_{\phi_m}\left(l''\left(h,X,Y\right)\right)\right| > \epsilon - a_m - c_m\right) \\
    & \leq \frac{\Var_\phi\left(\frac{\partial^2}{\partial^2\lambda}\ln\left(\frac{\partial^2}{\partial x\partial y}H\left(X_i,Y_i\right)_{\lambda + \frac{h}{\sqrt{k}}}\right)\right)}{\left(\epsilon - a_m - c_m\right)^2k}.
\end{align*}

Since the variance in the numerator is convergent and thus bounded due to the note following \textcolor{black}{Lemma}~\ref{Convergence of expectations}, this converges to zero as $k = k_n \rightarrow \infty$ and the theorem is proven.
\qed

\textcolor{black}{The above lemma assures us that a Taylor Series expansion of order 2 is sufficient to approximate the localized log-likelihood function. We can use this approximation to analyze the bias present in the likelihood function as it converges to that of a H\"usler-Reiss distribution.}

\begin{redlemma}\label{bias-variance theorem}
Assume the same conditions as \textcolor{black}{Lemma}~\ref{op theorem} and define
\begin{equation*}
    A:= \lim_{m\rightarrow \infty}\sqrt{k}\mathbb{E}_{\phi_m}\left(\frac{\partial}{\partial \lambda}\ln\left(h_{\lambda}\left(X,Y\right)\right)\right).
\end{equation*}
Then, $A$ is finite if and only if
\begin{equation*}
    \lim_{m\rightarrow\infty}\frac{\sqrt{k}}{b_m^2}, \hspace{5mm} \lim_{m\rightarrow\infty}\sqrt{k}\left(\lambda - b_m\sqrt{\frac{1-\rho_m}{1+\rho_m}}\right) 
\end{equation*}
are both finite. When $A$ is finite, we have
\begin{equation}\label{First Integral op}
    \frac{\partial}{\partial h}\tilde{L}\left(h\right) = \frac{\partial}{\partial h}\tilde{L}\left(0\right) - hI_{\lambda} + o_p\left(1\right) 
\end{equation}
\begin{equation}\label{Second Integral op}
    \tilde{L}\left(h\right) = \tilde{L}\left(0\right) + h\frac{\partial}{\partial h}\tilde{L}\left(0\right) - \frac{1}{2}h^2I_{\lambda} + o_p\left(1\right) 
\end{equation}
as $n \rightarrow \infty$. Here, $\frac{\partial}{\partial h}\tilde{L}\left(0\right) = \frac{1}{\sqrt{k}}\sum_{i=1}^k\frac{\partial}{\partial \lambda}\ln\left(h_{\lambda}\left(x_i,y_i\right)\right) \xrightarrow{d} N\left(A,I_\lambda\right)$.
\end{redlemma}

\noindent \textit{Proof.}
Equations~\ref{First Integral op} and~\ref{Second Integral op} are a direct result of integrating (\ref{Op}). What remains to show is the convergence of $\frac{\partial}{\partial h}\tilde{L}\left(0\right)$. Note from the form of $\frac{\partial}{\partial \lambda}\ln\left(h_{\lambda}\left(x_i,y_i\right)\right)$ (as made explicit in Section~\ref{Estimation Under the HR Model}) each term may once again be bounded by a term of the form $Ce^{-kx}e^{-jy}P\left(x-y\right)$, where $P\left(t\right)$ is some polynomial of $t$. Then the first and second moments of $\frac{\partial}{\partial \lambda}\ln\left(h_{\lambda}\left(x_i,y_i\right)\right)$ are convergent and bounded, which can be shown by the same method as used in \textcolor{black}{Lemma}~\ref{Convergence of expectations}. That is, the relevant variance and expectation are bounded, and we show that this implies that the Lindeberg condition for triangular arrays is satisfied (\cite{feller45}). Define $Var\left(\frac{\partial}{\partial\lambda}\ln\left(h_{\lambda}\left(X,Y\right)\right)\right) = \sigma^2_{k,i} = \sigma^2_{k}$, $s_k^2 = \sum_{i=1}^k\sigma^2_{k,i} = k\sigma^2_{k}$ and $\mu_k = \mathbb{E}_{\phi_m}\left(\frac{\partial}{\partial\lambda}\ln\left(h_{\lambda}\left(X,Y\right)\right)\right)$. In our case, the Lindeberg condition for triangular arrays reduces to

\begin{equation*}
    \lim_{n\rightarrow\infty}\frac{1}{\sigma_k^2}\mathbb{E}_{\phi_m}\left(\left(\frac{\partial}{\partial\lambda}\ln\left(h_{\lambda}\left(X,Y\right)\right) - \mu_k\right)^21_{\left|\frac{\partial}{\partial\lambda}\ln\left(h_{\lambda}\left(X,Y\right)\right) - \mu_k\right| > \epsilon s_k^2}\right) = 0.
\end{equation*}
Since $s_k^2 \rightarrow \infty$ as $k \rightarrow \infty$, we know that the pointwise limit of this function must be $0$. We also know that, as stated above, we may use the same DCT argument in \textcolor{black}{Lemma}~\ref{Convergence of expectations} to show that this expectation will converge to that pointwise limit, and thus the condition is satisfied. This gives us:
\begin{equation*}
    \frac{1}{\sqrt{k}}\sum_{i=1}^k\frac{\partial}{\partial \lambda}\ln\left(h_{\lambda}\left(x_i,y_i\right)\right) - \sqrt{k}\mathbb{E}_{\phi_m}\left(\frac{\partial}{\partial \lambda}\ln\left(\frac{\partial^2}{\partial x\partial y}H_{\lambda}\left(X,Y\right)\right)\right) \xrightarrow{D} N\left(0,I_\lambda\right)
\end{equation*}

Thus, our only remaining point is to show the conditions for $A$ to be finite. We first notice that $\frac{\partial}{\partial \lambda}\ln\left(\frac{\partial^2}{\partial x\partial y}H_{\lambda}\left(x,y\right)\right)$ is continuous in $\lambda$ for all $\lambda \in \left(0, \infty\right)$. By the Leibniz Integral Rule, we may pass the partial derivative with respect to $\lambda$ inside the expectation, giving us $\mathbb{E}_{\lambda}\left(\frac{\partial}{\partial \lambda}\ln\left(\frac{\partial^2}{\partial x\partial y}H_{\lambda}\left(X,Y\right)\right)\right) = 0$. As previously stated, we know that $\mathbb{E}_{\phi_m}\left(\frac{\partial}{\partial \lambda}\ln\left(\frac{\partial^2}{\partial x\partial y}H_{\lambda}\left(X,Y\right)\right)\right) \rightarrow \mathbb{E}_{\lambda}\left(\frac{\partial}{\partial \lambda}\ln\left(\frac{\partial^2}{\partial x\partial y}H_{\lambda}\left(X,Y\right)\right)\right) = 0$ by DCT, thus we need only to find the balance between the growth rates of $k$ and $m$.

We begin by rewriting the H\"usler-Reiss distribution as

\begin{equation}
    H_\lambda\left(x,y\right) = \exp\left\{-e^{-x}-\displaystyle\int_{y}^{\infty}\Phi\left(\lambda + \frac{x-z}{2\lambda}\right)e^{-z}dz\right\},
\end{equation}
and the distribution of the maxima of the rescaled normals as

\begin{align*}
    \Phi_{\rho_m}\left(u_m\left(x\right),u_m\left(y\right)\right) =& \exp\left\{-\frac{1}{m}\left[\int_y^{\infty}\Phi\left(\textcolor{black}{q_m(x,z)}\right)e^{-z-\frac{z^2}{2b_m^2}}dz + \right. \int_{x}^{\infty}e^{-t-\frac{t^2}{2b_m^2}}dt\right]\\
    &\left.-\sum_{i=2}^{\infty}\frac{\left(1-\Phi_{\rho_m}\left(u_m\left(x\right),u_m\left(y\right)\right)\right)^i}{i}\right\}.
\end{align*}

We write:
\begin{equation*}
\Biggr\rvert\mathbb{E}_{\phi_m}\left(\frac{\partial}{\partial \lambda}\ln\left(h_{\lambda}\left(X,Y\right)\right)\right)-0\Biggr\rvert=\Biggr\rvert\mathbb{E}_{\phi_m}\left(\frac{\partial}{\partial \lambda}\ln\left(h_{\lambda}\left(X,Y\right)\right)\right) - \mathbb{E}_{\lambda}\left(\frac{\partial}{\partial \lambda}\ln\left(h_{\lambda}\left(X,Y\right)\right)\right)\Biggr\rvert.
\end{equation*}
We further split this into the two terms:

\begin{align}\label{Split Expectation 1}
    \iint_{\mathbb{R}^2}\frac{\partial}{\partial \lambda}\ln\left(h_\lambda\left(x,y\right)\right)&\left[\Phi_{\rho_m}\left(u_m\left(x\right),u_m\left(y\right)\right)^{m-1}e^{-x-\frac{1}{2}\left(\frac{x}{b_m}\right)^2}\phi\left(\textcolor{black}{q_m(y,x)}\right)\frac{1}{b_m\sqrt{1-\rho_m^2}}\right.\\
    -&\left. H_\lambda\left(x,y\right)e^{-x}\phi\left(\lambda + \frac{y-x}{2\lambda}\right)\frac{1}{2\lambda}\right]dxdy\nonumber
\end{align}
and
\begin{align}\label{Split Expectation 2}
    &\iint_{\mathbb{R}^2}\frac{\partial}{\partial \lambda}\ln\left(h_\lambda\left(x,y\right)\right)\nonumber\\
    &\left[\Phi_{\rho_m}\left(u_m\left(x\right),u_m\left(y\right)\right)^{m-2}\frac{m-1}{m}e^{-x-y-\frac{1}{2}\left(\frac{x}{b_m}\right)^2-\frac{1}{2}\left(\frac{y}{b_m}\right)^2}\Phi\left(\textcolor{black}{q_m(y,x)}\right)\Phi\left(\textcolor{black}{q_m(x,y)}\right)\right.\nonumber\\
    &-\left. H_\lambda\left(x,y\right)e^{-x-y}\Phi\left(\lambda + \frac{y-x}{2\lambda}\right)\Phi\left(\lambda + \frac{x-y}{2\lambda}\right)\right]dxdy.
\end{align}
We begin with (\ref{Split Expectation 1}). We write
\begin{equation*}
    e^{-\frac{1}{2}\left(\frac{t}{b_m}\right)^2} = 1 - \frac{t^2}{2b_m^2} + e^{-\nu_t}\frac{t^4}{4b_m^4},
\end{equation*}
where $\frac{t^2}{2b_m^2} > \nu_t > 0$ is the value such that the final term is the Lagrange Remainder of order 2. Then we rewrite

\begin{align*}
    &\Phi_{\rho_m}\left(u_m\left(x\right),u_m\left(y\right)\right)^{m-1} = \exp\left\{-\frac{m-1}{m}\left[\int_y^{\infty}\Phi\left(\textcolor{black}{q_m(x,z)}\right)e^{-z-\frac{z^2}{2b_m^2}}dz \right.\right.\\
    & \hspace{10mm}+\left.\int_{x}^{\infty}e^{-t-\frac{t^2}{2b_m^2}}dz\right]\left.-(m-1)\sum_{i=2}^{\infty}\frac{\left(1-\Phi_{\rho_m}\left(u_m\left(x\right),u_m\left(y\right)\right)\right)^i}{i}\right\}\\
    &=\exp\left\{-\frac{m-1}{m}\left[\int_y^{\infty}\Phi\left(\textcolor{black}{q_m(x,z)}\right)e^{-z}\left(1 - \frac{z^2}{2b_m^2} + e^{-\nu_z}\frac{z^4}{4b_m^4}\right)dz + \right.\right.\\
    & \hspace{10mm}\left.\left. e^{-x}\left(1 - \frac{x^2+2x+2}{2b_m^2}\right) + \frac{1}{4b_m^4}\int_{x}^{\infty}e^{-t-\nu_t}t^4dt\right]-\left(m-1\right)\sum_{i=2}^{\infty}\frac{\left(1-\Phi_{\rho_m}\left(u_m\left(x\right),u_m\left(y\right)\right)\right)^i}{i}\right\}.
\end{align*}
We introduce the notation

\begin{equation}
    \Phi\left(a\right) = \Phi\left(b\right) + \phi\left(b\right)\left(b-a\right) - \frac{1}{2}\gamma_{a,b}\phi\left(\gamma_{a,b}\right)\left(b-a\right)^2
\end{equation}
where $\min\left\{a,b\right\}<\gamma_{a,b}<\max\left\{a,b\right\}$ make the final term the Lagrange Remainder of the expansion of $\Phi\left(a\right)$ about $b$ of order 2. Going forward, we drop the subscripts from $\gamma$ for the sake of notation. Using this formula to continue our calculations we see

\begin{align*}
    &\int_y^{\infty}\Phi\left(\textcolor{black}{q_m(x,z)}\right)e^{-z}\left(1 - \frac{z^2}{2b_m^2} + e^{-\nu_z}\frac{z^4}{4b_m^4}\right)dz\\
    &= \int_y^{\infty}\left(\Phi\left(\lambda + \frac{x-z}{2\lambda}\right) + \phi\left(\lambda + \frac{x-z}{2\lambda}\right)\left(\textcolor{black}{q_m(x,z)} - \lambda - \frac{x-z}{2\lambda}\right)\right.\\
    & \hspace{10mm}\left.-\frac{1}{2}\gamma_z\phi\left(\gamma_z\right)\left(\textcolor{black}{q_m(x,z)} - \lambda - \frac{x-z}{2\lambda}\right)^2\right)e^{-z}\left(1 - \frac{z^2}{2b_m^2} + e^{-\nu_z}\frac{z^4}{4b_m^4}\right)dz\\
    &= \int_y^{\infty}\Phi\left(\lambda + \frac{x-z}{2\lambda}\right)e^{-z}dz  + \frac{1}{4b_m^4}\int_y^{\infty}\Phi\left(\lambda + \frac{x-z}{2\lambda}\right)e^{-z}e^{-\nu_z}z^4dz\\
    & \hspace{10mm}- \frac{1}{2b_m^2}\int_y^{\infty}\Phi\left(\lambda + \frac{x-z}{2\lambda}\right)e^{-z}z^2dz+ \int_y^\infty\left(\phi\left(\lambda + \frac{x-z}{2\lambda}\right)\left(\textcolor{black}{q_m(x,z)} - \lambda - \frac{x-z}{2\lambda}\right)\right.\\
    & \hspace{10mm}\left.-\frac{1}{2}\gamma_z\phi\left(\gamma_z\right)\left(\textcolor{black}{q_m(x,z)} - \lambda - \frac{x-z}{2\lambda}\right)^2\right)e^{-z}\left(1 - \frac{z^2}{2b_m^2} + e^{-\nu_z}\frac{z^4}{4b_m^4}\right)dz.
\end{align*}
This leads us to the result 

\begin{align*}
    &\frac{\Phi_{\rho_m}\left(u_m\left(x\right),u_m\left(y\right)\right)^{m-1}}{H_\lambda\left(x,y\right)} = \exp\left\{f\left(x,y,m\right)\right\} :=\\
    & \exp\left\{-\frac{m-1}{m}\left[-\frac{1}{2b_m^2}\int_y^{\infty}\Phi\left(\lambda + \frac{x-z}{2\lambda}\right)e^{-z}z^2dz + \frac{1}{4b_m^4}\int_y^{\infty}\Phi\left(\lambda + \frac{x-z}{2\lambda}\right)e^{-z}e^{-\nu_z}z^4dz\right.\right.\\
    & \hspace{10mm}+ \int_y^\infty\left(\phi\left(\lambda + \frac{x-z}{2\lambda}\right)\left(\textcolor{black}{q_m(x,z)} - \lambda - \frac{x-z}{2\lambda}\right)\right.\\
    & \hspace{10mm}\left.-\frac{1}{2}\gamma_z\phi\left(\gamma_z\right)\left(\textcolor{black}{q_m(x,z)} - \lambda - \frac{x-z}{2\lambda}\right)^2\right)e^{-z}\left(1 - \frac{z^2}{2b_m^2} + e^{-\nu_z}\frac{z^4}{4b_m^4}\right)dz \\
    & \hspace{10mm}\left.\left.- e^{-x}\frac{x^2+2x+2}{2b_m^2}+ \frac{1}{4b_m^4}\int_{x}^{\infty}e^{-t-\nu_t}t^4dt\right]+\frac{1}{m}\left[e^{-x}+\int_{y}^{\infty}\Phi\left(\lambda + \frac{x-z}{2\lambda}\right)e^{-z}dz\right]\right.\\
    & \hspace{10mm}\left.-\left(m-1\right)\sum_{i=2}^{\infty}\frac{\left(1-\Phi_{\rho_m}\left(u_m\left(x\right),u_m\left(y\right)\right)\right)^i}{i}\right\}.
\end{align*}

\noindent Returning to (\ref{Split Expectation 1}),

\begin{align*}
    &\iint_{\mathbb{R}^2}\frac{\partial}{\partial \lambda}\ln\left(h_\lambda\left(x,y\right)\right)H_\lambda\left(x,y\right)\phi\left(\lambda + \frac{y-x}{2\lambda}\right)\frac{1}{2\lambda}e^{-x}\\
    &\hspace{10mm}\left[\frac{\Phi_{\rho_m}\left(u_m\left(x\right),u_m\left(y\right)\right)^{m-1}\phi\left(\textcolor{black}{q_m(y,x)}\right)2\lambda e^\frac{-x^2}{2b_m^2}}{H_\lambda\left(x,y\right)\phi\left(\lambda + \frac{y-x}{2\lambda}\right)b_m\sqrt{1-\rho_m^2}}-1\right]dxdy\\
    &=\iint_{\mathbb{R}^2}\frac{\partial}{\partial \lambda}\ln\left(h_\lambda\left(x,y\right)\right)H_\lambda\left(x,y\right)\phi\left(\lambda + \frac{y-x}{2\lambda}\right)\frac{1}{2\lambda}e^{-x}\\
    &\hspace{10mm}\left[\exp\left\{f\left(x,y,m\right)-\frac{x^2}{2b_m^2} - \frac{1}{2}\left(\textcolor{black}{q_m(y,x)}\right)^2 + \frac{1}{2}\left(\lambda + \frac{y-x}{2\lambda}\right)^2- \ln\left(\frac{b_m\sqrt{1-\rho_m^2}}{2\lambda}\right)\right\}-1\right]dxdy\\
    &=\iint_{\mathbb{R}^2}\frac{\partial}{\partial \lambda}\ln\left(h_\lambda\left(x,y\right)\right)H_\lambda\left(x,y\right)\phi\left(\lambda + \frac{y-x}{2\lambda}\right)\frac{1}{2\lambda}e^{-x}\left[f\left(x,y,m\right) - \frac{1}{2}\left(\textcolor{black}{q_m(y,x)}\right)^2\right.\\
    &\hspace{10mm}+ \frac{1}{2}\left(\lambda + \frac{y-x}{2\lambda}\right)^2 - \ln\left(\frac{b_m\sqrt{1-\rho_m^2}}{2\lambda}\right) - \frac{x^2}{2b_m^2} + e^{\zeta_a}\left(f\left(x,y,m\right) - \frac{1}{2}\left(\textcolor{black}{q_m(y,x)}\right)^2\right.\\
    &\left.\left.\hspace{10mm}+ \frac{1}{2}\left(\lambda + \frac{y-x}{2\lambda}\right)^2 - \ln\left(\frac{b_m\sqrt{1-\rho_m^2}}{2\lambda}\right)-\frac{x^2}{b_m^2}\right)^2\right]dxdy,
\end{align*}
where $\zeta_a$ is the appropriate Lagrange remainder for the expansion of the exponential term about 0. To show the convergence rate, we need only to assure that each term in the integral is integrable and converges to zero. We will show the absolute value of each term is integrable to simplify calculations. In the body of this proof, we will show this only for terms which may determine the rate of convergence, relegating any terms which necessarily converge faster to the appendix. This means that the terms accounted for here are:

\begin{itemize}
    \item [Ia] $= \displaystyle\frac{1}{2b_m^2}\int_y^{\infty}\Phi\left(\lambda + \frac{x-z}{2\lambda}\right)e^{-z}z^2dz$,
    \item [IIa] $= -\displaystyle\int_y^{\infty}\phi\left(\lambda + \frac{x-z}{2\lambda}\right)\left(\textcolor{black}{q_m(x,z)} - \lambda - \frac{x-z}{2\lambda}\right)e^{-z}dz$,
    \item [IIIa] $= \displaystyle e^{-x}\frac{x^2+2x+2}{2b_m^2}$,
    \item [IVa] $= -\displaystyle\ln\left(\frac{b_m\sqrt{1-\rho_m^2}}{2\lambda}\right)$,
    \item [Va] $= -\displaystyle\frac{x^2}{2b_m^2}$,
    \item [VIa] $\displaystyle= -\frac{1}{2}\left[\textcolor{black}{q_m(x,y)}^2 - \left(\lambda - \frac{x-y}{2\lambda}\right)^2\right]$.
\end{itemize}

To show that the above 6 terms are finite, we argue that each term is of the form $Ce^{-ax-by}x^{k_1}y^{k_2}$. We once again notice that each term in $\frac{\partial}{\partial \lambda}\ln\left(h_\lambda\left(x,y\right)\right)$ can be bounded by a term of the form $Ce^{-ax-by}\phi\left(\lambda + \frac{y-x}{2\lambda}\right)P\left(x-y\right)$ where $a,b\geq0$ and $P$ is a polynomial. These two facts combined will mean that each term above will be integrable by Lemma~\ref{integral}. Terms IIIa, Va and VIa are already in this form, and term IVa is a constant with respect to our integral, and thus does not affect the integrability of the integrand.

To handle term Ia, we use the bounds:

\begin{equation*}
    \text{Ia} = \frac{1}{2b_m^2}\int_y^{\infty}\Phi\left(\lambda + \frac{x-z}{2\lambda}\right)e^{-z}z^2dz \leq \frac{1}{2b_m^2}\int_y^\infty e^{-z}z^2dz = e^{-y}\frac{y^2 + 2y + 2}{2b_m^2}.
\end{equation*}

To analyze term IIa we first note

\begin{equation*}
    \textcolor{black}{q_m(x,y)}- \lambda - \frac{\textcolor{black}{x-y}}{2\lambda} = b_m\sqrt{\frac{1-\rho_m}{1+\rho_m}} - \lambda + \left(x-y\right)\left[\frac{1}{b_m\sqrt{\left(1-\rho_m^2\right)}} - \frac{1}{2\lambda}\right] + x \sqrt{\frac{1-\rho_m}{1+\rho}}b_m^{-1}.
\end{equation*}
We notice that the first and second pairs of terms need not converge at the rate $\frac{1}{log\left(m\right)}$ or faster as all other terms have so far. Thus if this term converges slower than $1/b_m^2$, this will dictate our asymptotic bias. We continue, showing that this term is indeed integrable:

\begin{align*}
    \text{IIa} &= \int_y^{\infty}\phi\left(\lambda + \frac{x-z}{2\lambda}\right)\left(\textcolor{black}{q_m(x,z)} - \lambda - \frac{x-z}{2\lambda}\right)e^{-z}dz\\
    &\leq \int_y^{\infty}e^{-z}\left(b_m\sqrt{\frac{1-\rho_m}{1+\rho_m}} - \lambda + \left(x-z\right)\left[\frac{1}{b_m\sqrt{\left(1-\rho_m^2\right)}} - \frac{1}{2\lambda}\right] + z \sqrt{\frac{1-\rho_m}{1+\rho}}b_m^{-1}\right)dz\\
    &= e^{-y}\left(b_m\sqrt{\frac{1-\rho_m}{1+\rho_m}} - \lambda + \left(x-y-1\right)\left[\frac{1}{b_m\sqrt{\left(1-\rho_m^2\right)}} - \frac{1}{2\lambda}\right] + (y+1) \sqrt{\frac{1-\rho_m}{1+\rho}}b_m^{-1}\right),
\end{align*}
which again allows us to utilize the form of Lemma~\ref{integral}.

The analysis of (\ref{Split Expectation 2}) is generally similar, but contains a few unique details. First, we may ignore the term $\frac{m-1}{m} = 1- \frac{1}{m}$, as factoring out $\frac{1}{m}$ is already a faster rate of convergence than other terms. So, understanding that this term doesn't affect our analysis, we drop it in further analysis. Adopting the notation that $f'\left(x,y,m\right)$ is the same as $f\left(x,y,m\right)$ but with the quantity $m-2$ in place of $m-1$, our equation becomes:

\begin{align*}
    &\iint_{\mathbb{R}^2}\frac{\partial}{\partial \lambda}\ln\left(h_\lambda\left(x,y\right)\right)H_\lambda\left(x,y\right)e^{-x-y}\left[-\Phi\left(\lambda + \frac{y-x}{2\lambda}\right)\Phi\left(\lambda + \frac{x-y}{2\lambda}\right)\right.\\
    &\hspace{10mm}+\left. e^{f'\left(x,y,m\right)-\frac{1}{2}\left(\frac{x}{b_m}\right)^2-\frac{1}{2}\left(\frac{y}{b_m}\right)^2}\Phi\left(\textcolor{black}{q_m(y,x)}\right)\Phi\left(\textcolor{black}{q_m(x,y)}\right)\right]dxdy\\
    &=\iint_{\mathbb{R}^2}\frac{\partial}{\partial \lambda}\ln\left(h_\lambda\left(x,y\right)\right)H_\lambda\left(x,y\right)e^{-x-y}\left[-\Phi\left(\lambda + \frac{y-x}{2\lambda}\right)\Phi\left(\lambda + \frac{x-y}{2\lambda}\right)\right.\\
    &\hspace{10mm}+ \left(1 + f'\left(x,y,m\right) - \frac{1}{2b_m^2}\left(x^2 + y^2\right) + \frac{1}{2}\left(f'\left(x,y,m\right) - \frac{1}{2b_m^2}\left(x^2 + y^2\right)\right)^2e^{\zeta_b}\right)\\
    &\hspace{10mm}\left(\Phi\left(\lambda + \frac{y-x}{2\lambda}\right) + \left(\textcolor{black}{q_m(y,x)} - \lambda - \frac{y-x}{2\lambda}\right)\phi\left(\lambda + \frac{y-x}{2\lambda}\right)\right.\\
    &\left.\hspace{10mm}- \frac{1}{2}\left(\textcolor{black}{q_m(y,x)} - \lambda - \frac{y-x}{2\lambda}\right)^2\gamma_1\phi\left(\gamma_1\right)\right)\\
    &\hspace{10mm}\left(\Phi\left(\lambda + \frac{x-y}{2\lambda}\right) + \left(\textcolor{black}{q_m(x,y)} - \lambda - \frac{x-y}{2\lambda}\right)\phi\left(\lambda + \frac{x-y}{2\lambda}\right)\right.\\
    &\hspace{10mm}-\left.\left. \frac{1}{2}\left(\textcolor{black}{q_m(x,y)} - \lambda - \frac{x-y}{2\lambda}\right)^2\gamma_2\phi\left(\gamma_2\right)\right)\right]dxdy,
\end{align*}
where $\min\left(\textcolor{black}{q_m(y,x)},\lambda + \frac{y-x}{2\lambda}\right) \leq \gamma_1 \leq \max\left(\textcolor{black}{q_m(y,x)},\lambda + \frac{y-x}{2\lambda}\right)$ and $\gamma_2$ is defined similarly. 
Note that, once factored out, $\Phi\left(\lambda + \frac{y-x}{2\lambda}\right)\Phi\left(\lambda + \frac{x-y}{2\lambda}\right)$ will cancel.
Again, in the body of this proof we only show the integrability of terms which may affect our convergence rate. The terms we consider are:

\begin{itemize}
    \item [Ib] $= \displaystyle\Phi\left(\lambda + \frac{y-x}{2\lambda}\right)\Phi\left(\lambda + \frac{x-y}{2\lambda}\right)\frac{1}{2b_m^2}\int_y^{\infty}\Phi\left(\lambda + \frac{x-z}{2\lambda}\right)e^{-z}z^2dz$,
    \item [IIb] $= -\displaystyle\Phi\left(\lambda + \frac{x-y}{2\lambda}\right)\Phi\left(\lambda + \frac{y-x}{2\lambda}\right)\int_y^{\infty}\phi\left(\lambda + \frac{x-z}{2\lambda}\right)\left(\textcolor{black}{q_m(x,z)} - \lambda - \frac{x-z}{2\lambda}\right)e^{-z}dz$,
    \item [IIIb] $= \displaystyle\Phi\left(\lambda + \frac{x-y}{2\lambda}\right)\Phi\left(\lambda + \frac{y-x}{2\lambda}\right)e^{-x}\frac{x^2+2x+2}{2b_m^2}$,
    \item [IVb] $= \displaystyle\Phi\left(\lambda + \frac{x-y}{2\lambda}\right)\Phi\left(\lambda + \frac{y-x}{2\lambda}\right)\frac{-x^2}{2b_m^2}$,
    \item [Vb] $= \displaystyle\Phi\left(\lambda + \frac{x-y}{2\lambda}\right)\Phi\left(\lambda + \frac{y-x}{2\lambda}\right)\frac{-y^2}{2b_m^2}$,
    \item [VIb] $= \displaystyle\Phi\left(\lambda + \frac{x-y}{2\lambda}\right)\left(\textcolor{black}{q_m(y,x)} - \lambda - \frac{y-x}{2\lambda}\right)\phi\left(\lambda + \frac{y-x}{2\lambda}\right)$,
    \item [VIIb] $= \displaystyle \Phi\left(\lambda + \frac{y-x}{2\lambda}\right)\left(\textcolor{black}{q_m(x,y)} - \lambda - \frac{x-y}{2\lambda}\right)\phi\left(\lambda + \frac{x-y}{2\lambda}\right)$.
\end{itemize}

Terms Ib through Vb are the same terms (or the same term with $y$ in place of $x$) as in the previous analysis, with a leading $\Phi\left(\lambda + \frac{y-x}{2\lambda}\right)\Phi\left(\lambda + \frac{x-y}{2\lambda}\right)$. Each of these terms may be written as $\Phi\left(\lambda + \frac{x-y}{2\lambda}\right)\Phi\left(\lambda + \frac{y-x}{2\lambda}\right)e^{-ax-by}x^{k_1}y^{k_2}$. Lemma~\ref{appendix4} shows that each of these terms is finite. Similarly, VIb and VIIb already contain the term $\phi(\lambda + \frac{x-y}{2\lambda})$, so we can use the same argument as in Lemma~\ref{integral} and thus the theorem is proven. We note that this theorem also gives the form of $A$, which is made explicit in the appendix.
\qed

Finally, we come to the proof of Theorem~\ref{MLE Form}. The proof of this theorem is nearly identical to the proof of Theorem 2.2 in \textcites{dombry19}, the only difference being the parameter space in \textcites{dombry19} is $\mathbb{R}^3$, while ours is $\mathbb{R}$. What follows is a simplified version of their Theorem 2.2, reduced to the 1-dimensional case.

\noindent \textit{Proof \textcolor{black}{of Theorem \ref{MLE Form}}.}
Consider the processes

\begin{equation}
    M_n\left(h\right) = \tilde{L}\left(h\right) - \tilde{L}\left(0\right)
\end{equation}
and 
\begin{equation}
    M\left(h\right) = h\left(A + G\right) - \frac{1}{2}h^2I_\lambda
\end{equation}
where $G$ is a 0-mean Gaussian random variable with covariance $I_\lambda$. Let $H_n$ be a closed ball in $\mathbb{R}$ centered at 0 and with radius $r_k \rightarrow \infty$ such that $r_k = o\left(k^{-\frac{1}{2}}\right)$. Define $\hat{h}_n = \text{argmax} M_n\left(h\right)$ for $h \in H_n$. By \textcolor{black}{Lemma}~\ref{bias-variance theorem}, for any compact $K \subset \mathbb{R}$, $M_n$ converges in distribution to $M$ in $l^{\infty}\left(K\right)$ as $n \rightarrow \infty$. $M$ is continuous and has the unique maximizer

\begin{equation}
    \hat{h} = \text{argmax} M\left(h\right) = I_\lambda^{-1}\left(A+G\right).
\end{equation}

By a corollary of the Argmax theorem (\cite{vaart_1998}, Corollary 5.58), if $\hat{h}_n$ is tight, then it converges in distribution to $\hat{h}$. Thus all that remains is to show that $\hat{h}_n$ is tight. Fix $\epsilon >0$. Note that the form of $\hat{h}$ assures us that there must exist some $R > 0$ such that $\mathbb{P}\left(|\hat{h}| < R\right) > 1- \epsilon$. The form of $M$ gives us

\begin{equation*}
    M\left(h\right) = M\left(\hat{h}\right) - \frac{1}{2}\left(h-\hat{h}\right)^2I_\lambda,
\end{equation*}
and thus 

\begin{equation*}
    \max_{|\hat{h}-h|\geq 1}M\left(h\right) = M\left(\hat{h}\right) -\frac{1}{2}I_\lambda.
\end{equation*}
For $|\hat{h}| < R$

\begin{equation*}
    M\left(\hat{h}\right) - \max_{|h| = R + 1}M\left(h\right) = \max_{|h| \leq R} M\left(h\right) - \max_{|h| = R + 1}M(h) \geq \frac{1}{2} I_\lambda,
\end{equation*}
which, based on the selection of $R$, occurs with probability at least $1-\epsilon$. Define $K = \left\{h:|h| \leq R + 1\right\}$. Since $M_n$ converges weakly to $M$ in $l^{\infty}\left(K\right)$, for
sufficiently large $n$,

\begin{equation*}
	\max_{|h|\leq R}M_n\left(h\right) - \max_{|h|= R + 1}M_n\left(h\right) \geq .25I_\lambda
\end{equation*}
with probability at least $1 - 2\epsilon$. For sufficiently large $n$, $K \subset H_n$. Hence $M_n$
is strictly concave on $H_n$ by \textcolor{black}{Lemma}~\ref{Convergence of expectations} with probability
at least $1 - \epsilon$. Then, from the equation above, and the fact that $M_n$ is strictly concave, we
see that $\mathbb{P}\left(|\hat{h}_n| \leq R\right) \geq 1-3\epsilon$. Thus, $\hat{h}$ is tight. Also note,
that $\hat{h}_n$ is on the interior of $H_n$, thus

\begin{equation*}
	\frac{\partial}{\partial h}M_n\left(\hat{h}_n\right) = 0 \implies \frac{\partial}{\partial h}\tilde{L}\left(\hat{h}_n\right) = 0.
\end{equation*}

Define $\hat{\lambda}_n = \lambda + \frac{\hat{h}_n}{\sqrt{k}}$. Then if $\hat{h}_n$ maximizes the localized log-likelihood $\tilde{L}$, $\hat{\lambda}$ must maximize the log-likelihood $L$. Thus, with probability tending to one, $\hat{\lambda}$ is an MLE.
\qed

\section{Simulation Study}
\begin{figure}
\centering
\begin{minipage}{.45\textwidth}
\centering
     \includegraphics[width=.9\textwidth]{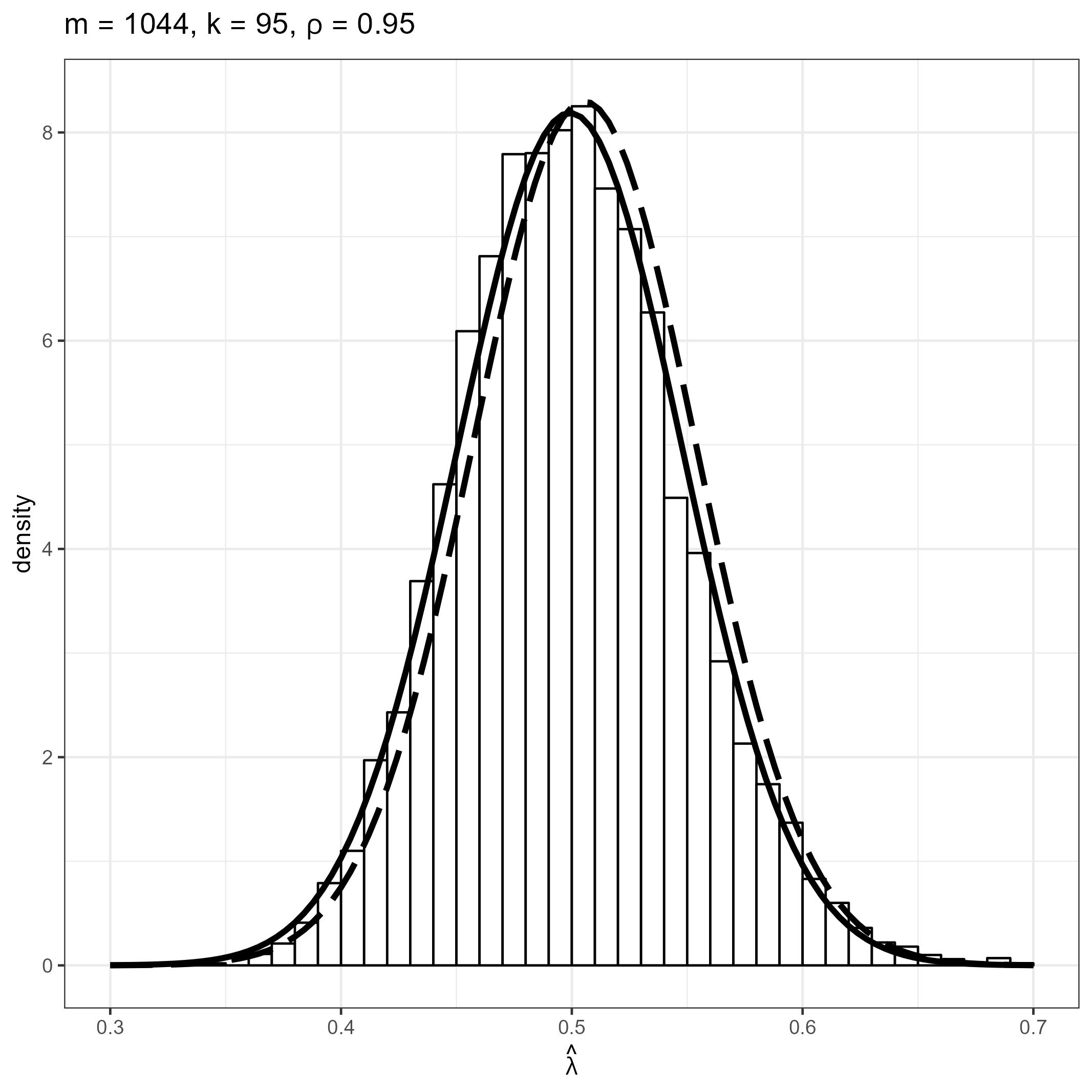}
\end{minipage}\hfill
\begin{minipage}{.45\textwidth}
    \centering
    \includegraphics[width=.9\textwidth]{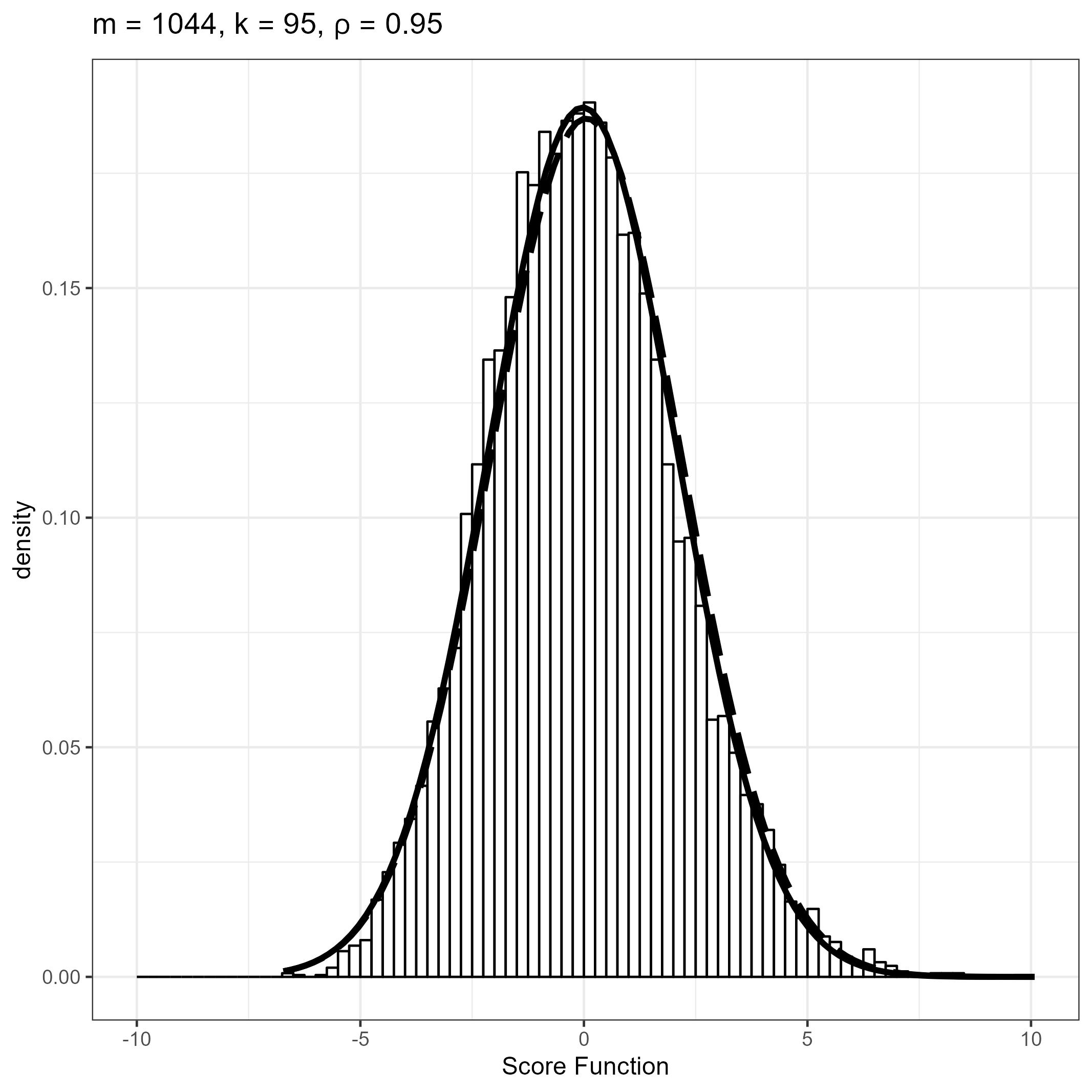}
\end{minipage}
\caption{Left: 10,000 replications of simulated $\hat{\lambda}$ with $n=\textcolor{black}{99,180}$. Right: Simulated values of the score function from the same sample. Dashed lines represent theoretical distributions, solid lines represent fitted distributions}
\label{tab:sim1}
\end{figure}

\begin{figure}
\centering
\begin{minipage}{.45\textwidth}
\centering
     \includegraphics[width=.9\textwidth]{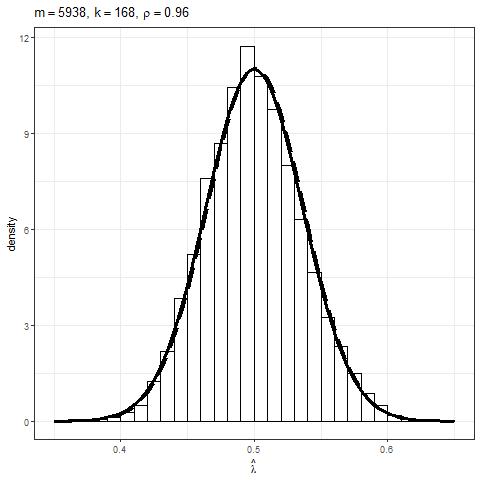}
\end{minipage}\hfill
\begin{minipage}{.45\textwidth}
    \centering
    \includegraphics[width=.9\textwidth]{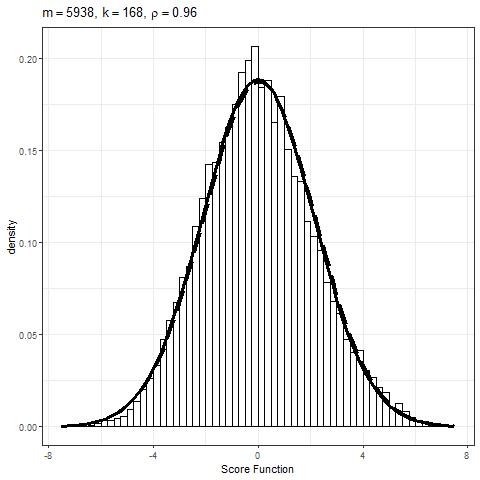}
\end{minipage}
\caption{Left: 10,000 replications of simulated $\hat{\lambda}$ with $n=\textcolor{black}{997,584}$. Right: Simulated values of the score function from the same sample. Dashed lines represent theoretical distributions, solid lines represent fitted distributions. \textcolor{black}{We note that the agreement between the fitted values is close enough that the lines are barely distinguishable; See Table \ref{tab:sim1} for the fitted values.}}
\label{sim2}
\end{figure}

\begin{figure}
\centering
\begin{minipage}{.45\textwidth}
\centering
     \includegraphics[width=.9\textwidth]{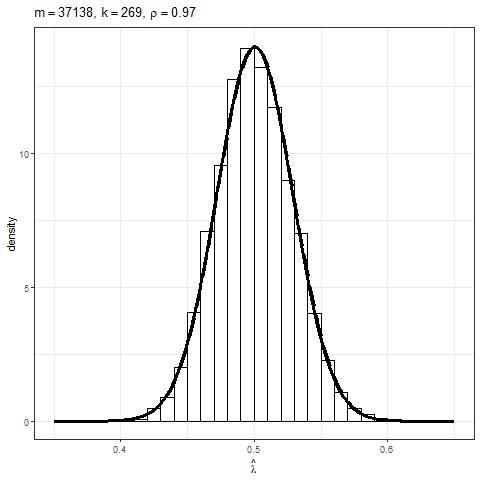}
\end{minipage}\hfill
\begin{minipage}{.45\textwidth}
    \centering
    \includegraphics[width=.9\textwidth]{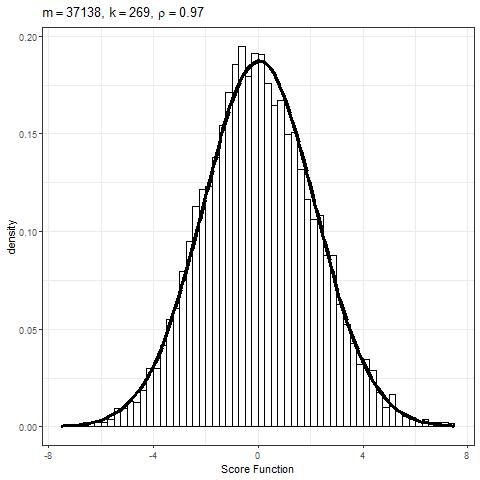}
\end{minipage}
\caption{ Left: 10,000 replications of simulated $\hat{\lambda}$ with $n=\textcolor{black}{9,990,122}$. Right: Simulated values of the score function from the same sample. Dashed lines represent theoretical distributions, solid lines represent fitted distributions. \textcolor{black}{We note that the agreement between the fitted values is close enough that the lines are barely distinguishable; See Table \ref{tab:sim1} for the fitted values.}}
\label{sim3}
\end{figure}

To illustrate the practical implications of the results, we run a simulation study. A simulation of the MLE (labelled $\hat{\lambda}$) is repeated 10,000 times. As shown below, the convergence to the limiting results is very slow. We mitigate this by first considering a sample size \textcolor{black}{of roughly} $n = 100,000$. \textcolor{black}{We say roughly here because we choose to somewhat simplify the problem by ensuring that each block is of equal size, and thus we adjust the overall sample size such that it divides evenly into the blocks}. After selecting a \textcolor{black}{rough} sample size, to determine the values for $m$ and $k$ we must select a formula for $\rho_m$. For convenience, we select a formula for $\rho_m$ such that $\lambda_m \textcolor{black}{:= b_m\sqrt{\frac{1-\rho_m}{1+\rho_m}}} = \lambda$ at each step:
\begin{equation} \label{rhoform}
    \rho_m = \frac{b_m^2-\lambda^2}{b_m^2 + \lambda^2}.
\end{equation} 
This assures us that the expectation of the score function converges to zero at a rate of $1/b_m^2$, up to a constant. \textcolor{black}{We then choose $k = b_m^4$ to satisfy the finite limit conditions of Theorem \ref{MLE Form}}. We can then use the definition of $b_m$ and $mk = n$ to solve for $b_m$ and subsequently $m$ and $k$. With this method we arrive at the values $m = 1,044$ and $k = 95$. \textcolor{black}{Although in principle, there is no restriction that each block be of exactly the same size so long as they produce the same $\lambda$, it is unclear what the best way to deal with unequal blocks is. We choose to round $k$ down so each of our blocks is of the same size, avoiding the problem altogether. As such, $mk = 99,180$ not the full $n=100,000$, but, for the purposes of these simulations, this difference is negligible}. Selecting $\lambda = .5$, this gives us $\rho_m = 0.95$. \textcolor{black}{The data are generated as bivariate normal random variables, with the mean and variance given by
\begin{equation}
    \mu = (-b_m^2,-b_m^2)^T, \quad \Sigma = b_m^2\begin{bmatrix}
1 & \rho_m \\
\rho_m & 1 
\end{bmatrix},
\end{equation}with the form of $\rho_m$ specified in (\ref{rhoform}). This provides the same effect as generating the data as standard normal random variables with correlation $\rho_m$, and subsequently applying the scaling. We generate $m = 1,044$ samples from this distribution, $k = 95$ times. The maxima are selected for each set of $m$.} The results of this simulation are shown in Figure~\ref{tab:sim1}. We report both the simulated $\hat{\lambda}$s as well as the values of the score functions, whose distribution (as made explicit in \textcolor{black}{Lemma}~\ref{bias-variance theorem}) depends on the same parameters $A$ and $I_{\lambda}$. Table~\ref{simtab} numerically summarizes the results of this simulation. We see close agreement between theoretical and experimental values for the variances of both the score function and $\hat{\lambda}$. The empirical bias term, however seems not to initially agree with the theoretical values. To investigate this disagreement, we repeat this simulation with increased \textcolor{black}{(rough)} sample sizes of $n=1,000,000$ and $n=10,000,000$, the results of which are also summarized in Table~\ref{simtab}, Figure~\ref{sim2} and Figure~\ref{sim3}, respectively. \textcolor{black}{We repeat the same process to find $k$ and $m$, again ensuring that our blocks are of equal size, and end with final sample sizes of 997,584 and 9,990,122 respectively.} We see that it takes until the final experiment for a 95\% confidence interval to contain our theoretical value for both the bias of the score function and $\hat{\lambda}$s. We additionally see that with each step the confidence interval shifts right toward our predicted theoretical value. This implies that convergence toward the theoretical bias value takes large sample sizes. We also note that the variances of these distributions are quite large compared to the size of the bias terms, and as such noisy estimates of these means are to be expected. \textcolor{black}{The expected noise is supported} by repeating these simulations, which tend to produce fairly different values for the means in many runs.

\textcolor{black}{The variability we see present in the simulations} underpins the importance of finding the correct balance between $k$ and $m$. We see through the previous simulations that an increase in $k$ does reduce the variability, but only slightly, and that large samples of $m$ are necessary to find accurate estimates of the bias term. To that end, we run the same simulation, with $n = 1,000,000$ but redefine the relationship between $k$ and $m$ to be $k = \frac{1}{4}b_m^4$. This gives us the values $k = 56$ and $m = 17,732$ \textcolor{black}{for a final sample size of $n=992,992$}. This simulation shows the closest agreement with the theoretical values thus far. While the performance of this simulation can at least partially be attributed to random variation, it exemplifies the importance of striking the right balance in the bias variance tradeoff, as we are able to get similar or better performance to the simulation with \textcolor{black}{nearly ten times the sample size}.

\textcolor{black}{The problem of balancing $k$ and $m$ is still a fairly tricky one. The finite limit conditions of Theorem~\ref{MLE Form} give us the condition $k = O(b_m^{4})$, but $\left(\lambda - b_m\sqrt{\frac{1-\rho_m}{1+\rho_m}}\right)$ may grow arbitrarily slowly or quickly and thus it is even hard to assure the the asymptotic bias will be finite. There is a so-called second order H\"usler-Reiss condition, as per \textcite{liao} in which $\left(\lambda - b_m\sqrt{\frac{1-\rho_m}{1+\rho_m}}\right) = O(b_m^{-4})$, but even this does not mitigate the fact that the value of the asymptotic bias depends on the true value of $\lambda$. However, in practice the asymptotic bias term is so small for reasonable values of $\lambda$ (roughly $(0,1]$), that the asymptotic bias is insignificant when compared to the bias present at smaller sample sizes. It is a result dating back to \textcite{FT28} that Guassian variables require large sample sizes for the maxima to converge to the limiting extreme value distribution. As such the simulations above use sample sizes which are unreasonable for any practical application. We note, however, that the sample sizes presented in this paper are not needed for accurate prediction of the true value of $\lambda$. Since our results are second order asymptotics, these huge samples a required to isolate the asymptotic bias such that we can support the theoretical findings with simulation.}

\newpage

\begin{table}

\begin{tabular}{|c|c|c|c|c|c|c|}
\hline
& \begin{tabular}[c]{@{}l@{}}Sample \\ Size\end{tabular}& \begin{tabular}[c]{@{}l@{}}Theoretical \\ Mean\end{tabular} & \begin{tabular}[c]{@{}l@{}}Simulated \\ Mean\end{tabular} & 95\% C.I. & \begin{tabular}[c]{@{}l@{}}Theoretical \\ Variance\end{tabular} & \begin{tabular}[c]{@{}l@{}}Simulated \\ Variance\end{tabular} \\ \hline
\multirow{3}{*}{$\displaystyle k=b_m^4$} & $\textcolor{black}{99,180}$ & 0.0117 & -0.0081 & (-0.0173,0.0013) &0.2196  & 0.2255\\ 
\cline{2-7} & $\textcolor{black}{997,584}$ & 0.0117 & -0.0029 & (-0.0121,0.0063) &0.2196  & 0.2209\\
\cline{2-7} & $\textcolor{black}{9,990,122}$ & 0.0117 & 0.0035 & (-0.0057,0.0127) & 0.2196 & 0.2187 \\ \hline
$\displaystyle k=\frac{1}{4}b_m^4$ & \textcolor{black}{992,992} & 0.0058 & 0.0038 &(-0.0055,0.0131)& 0.2196 & 0.2258\\
\hline
\end{tabular}
\captionof{table}{Results from four simulations, each with 10,000 replications.}
\label{simtab}
\end{table}

\textbf{Statements and Declarations:} Jan Hannig's research was supported in part by the National Science Foundation under Grant No. DMS-1916115, 2113404, and 2210337. Hank Flury was also support in part by the National Science Foundation under Grant No. DMS-2210337. The authors have no competing interests to report. This work did not involve humans and/or animals and as such has no approval committee. All code, including generated data can be found at the public GitHub repository: \url{https://github.com/fluryh/Asymptotic_Theory_for_Estimation-_of_the_Husler-Reiss_Distribution_via_Block_Maxima_Method}.
% Bibliography
\printbibliography

@Article{         FT28,
  author	= {Fisher, R.A. and Tippett, L.H.C.},
  title		= {Limiting forms of the frequency distribution of the largest or smallest member of a sample},
  journal	= {Proc. Cambridge Phil. Soc.},
  volume 	= {24},
  pages         = {180--190},
  year		= {1928}
}

@Article{         tawn88,
  author	= {Tawn, J.A.},
  title		= {Bivariate extreme value theory: Models and estimation},
  journal	= {Biometrika},
  volume        = {75(3)},
  pages         = {397--415},
  year		= {1988}
}

@Article{         tawn90,
  author	= {Tawn, J.A.},
  title		= {Modelling multivariate extreme value distributions},
  journal	= {Biometrika},
  volume        = {77(2)},
  pages         = {245--253},
  year		= {1990}
}

@Article{         huslerreiss89,
  author	= {H\"usler, J. and Reiss, R.-D.},
  title		= {Maxima of normal random variables: between independence and complete
dependence},
  journal	= {Statistics and Probability Letters},
  volume	= {7},
  pages		= {283--286},
  year		= {1989}
}

@Article{         davisonpadoanribatet,
  author	= {Davison, A.C. and Padoan, S.A. and Ribatet, M.},
  title		= {Statistical modeling of spatial extremes (with discussion)},
  journal	= {Statistical Science},
  volume        = {27},
  pages         = {161--186},
  year		= {2012}
}

@Article{         dombry19,
  author	= {Dombry, C. and Ferreira, A.},
  title		= {Maximum likelihood estimators based on the block maxima method},
  journal	= {Bernoulli},
  volume	= {25},
  pages		= {1690--1723},
  year		= {2019}
}

@Article{         liao,
  author	= {X. Liao and Z. Peng},
  title		= {{Convergence rate of maxima of bivariate Guassian arrays to the H}\"usler-{Reiss distribtuion}},
  journal	= {Statistics and Its Interface},
  volume	= {7},
  pages		= {351-362},
  year		= {2014}
}

@book{vaart_1998, place={Cambridge}, series={Cambridge Series in Statistical and Probabilistic Mathematics}, title={Asymptotic Statistics}, DOI={10.1017/CBO9780511802256}, publisher={Cambridge University Press}, author={{van der Vaart}, A. W.}, year={1998}, collection={Cambridge Series in Statistical and Probabilistic Mathematics}}

@article{engelke2015,
 ISSN = {13697412, 14679868},
 URL = {http://www.jstor.org/stable/24774733},
 author = {Sebastian Engelke and Alexander Malinowski and Zakhar Kabluchko and Martin Schlather},
 journal = {Journal of the Royal Statistical Society. Series B (Statistical Methodology)},
 number = {1},
 pages = {239--265},
 publisher = {[Royal Statistical Society, Wiley]},
 title = {Estimation of {H}\"usler–{R}eiss distributions and {B}rown–{R}esnick processes},
 urldate = {2023-04-21},
 volume = {77},
 year = {2015}
}

@book{abramowitz+stegun,
  address = {New York},
  author = {Abramowitz, Milton and Stegun, Irene A.},
  edition = {ninth Dover printing, tenth GPO printing},
  publisher = {Dover},
  title = {Handbook of Mathematical Functions with Formulas, Graphs, and Mathematical Tables},
  year = 1964
}

@article{feller45,
author = {W. Feller},
title = {{The fundamental limit theorems in probability}},
volume = {51},
journal = {Bulletin of the American Mathematical Society},
number = {11},
publisher = {American Mathematical Society},
pages = {800 -- 832},
year = {1945},
}

@article{Asadi_2015,
   title={Extremes on river networks},
   volume={9},
   ISSN={1932-6157},
   url={http://dx.doi.org/10.1214/15-AOAS863},
   DOI={10.1214/15-aoas863},
   number={4},
   journal={The Annals of Applied Statistics},
   publisher={Institute of Mathematical Statistics},
   author={Asadi, Peiman and Davison, Anthony C. and Engelke, Sebastian},
   year={2015},
   month=dec }
\newpage
\appendix
\section{Technical Proofs}

\begin{lemma}\label{appendix1}
    Integrals of the form 
    \begin{multline*}
        \iint_{\mathbb{R}\textcolor{black}{^2}}x^{k_1}y^{k_2}e^{-ax-by}\phi\left(\frac{\left(x-y\right) -\mu}{\sigma}\right)H_\lambda\left(x,y\right)
        \left(m-1\right)\sum_{i=2}^{\infty}\frac{\left(1-\Phi_{\rho_m}\left(u_m\left(x\right),u_m\left(y\right)\right)\right)^i}{i}dxdy
    \end{multline*} are finite for all $m$ such that $\rho_m, a+b+1,k_1,k_2 > 0$.
\end{lemma}

\noindent \textit{Proof.}
We first establish the bound
\begin{equation*}
    1-\Phi(u_m(x)) = \int_{u_m(x)}^{\infty}\phi(t)dt
    = \frac{1}{b_m}\phi(b_m)\int_x^\infty e^{-s-\frac{s^2}{2b_m^2}}ds
    \leq \frac{1}{m}e^{-x}
\end{equation*}
With this bound, and the trivial bound $1-\Phi_{\rho_m}(u_m(x),u_m(y)) \leq 1 - \Phi(u_m(x)) + 1 - \Phi(u_m(y))$, we derive the following bound.
\begin{align*}
    &\left(m-1\right)\sum_{i=2}^{\infty}\frac{\left(1-\Phi_{\rho_m}\left(u_m\left(x\right),u_m\left(y\right)\right)\right)^i}{i} \\
    &\leq\left(m-1\right)\left(1-\Phi_{\rho_m}\left(u_m\left(x\right),u_m\left(y\right)\right)\right)^2\sum_{i=0}^{\infty}\left(1-\Phi_{\rho_m}\left(u_m\left(x\right),u_m\left(y\right)\right)\right)^i\\
    &= \left(m-1\right)\frac{\left(1-\Phi_{\rho_m}\left(u_m\left(x\right),u_m\left(y\right)\right)\right)^2}{\Phi_{\rho_m}\left(u_m\left(x\right),u_m\left(y\right)\right)}\\
    &\leq \frac{(m-1)}{m^2}\left(e^{-x}+e^{-y}\right)^2\frac{1}{\Phi_{\rho_m}\left(u_m\left(x\right),u_m\left(y\right)\right)}\\
    &\leq \left(e^{-x}+e^{-y}\right)^2\frac{1}{\Phi_{\rho_m}\left(u_m\left(x\right),u_m\left(y\right)\right)},
\end{align*}
For $x,y>0$, $\Phi_{\rho_m}\left(x,y\right) > \frac{1}{4}$. For $x,y$ such that $\min\left(x,y\right)<0$ and $\rho_m > 0$
\begin{align*}
    \Phi_{\rho_m}\left(x,y\right) &\geq \Phi_{\rho_m}\left(\min\left(x,y\right),\min\left(x,y\right)\right)\\
    &= \mathbb{P}(X < \min(x,y)|Y < \min(x,y))\mathbb{P}(Y<\min(x,y))\\
    &\geq \mathbb{P}(X < \min(x,y))\mathbb{P}(Y<\min(x,y))\\
    & = \Phi(\min(x,y))^2.
\end{align*}

With this, we may return to the integral in the statement of Lemma~\ref{appendix1}. We drop the term $(e^{-x}+e^{-y})^2$. The only affect of this term is to replace increase the coefficients $a$ or $b$ by 2, or each by 1. This has no affect on the technique going forward, and so the term is dropped. We write the integral as:  

\begin{align*}
    &\iint x^{k_1}y^{k_2}e^{-ax-by}\phi\left(\frac{\left(x-y\right) -\mu}{\sigma}\right)H_\lambda\left(x,y\right)\frac{1}{\Phi_{\rho_m}\left(u_m\left(x\right),u_m\left(y\right)\right)}dxdy\\
    &\leq 4\iint_{\mathbb{R}^2\setminus\left(-b_m^2,\infty\right)^2} x^{k_1}y^{k_2}e^{-ax-by}\phi\left(\frac{\left(x-y\right) -\mu}{\sigma}\right)H_\lambda\left(x,y\right)dxdy\\
    &+ \iint_{\left(-b_m^2,\infty\right)^2} x^{k_1}y^{k_2}e^{-ax-by}\phi\left(\frac{\left(x-y\right) -\mu}{\sigma}\right)H_\lambda\left(x,y\right)\frac{1}{\Phi\left(u_m\left(\min\left(x,y\right)\right)\right)^2}dxdy\\
    &\leq 4\iint_{\mathbb{R}^2\setminus\left(-b_m^2,\infty\right)^2} x^{k_1}y^{k_2}e^{-ax-by}\phi\left(\frac{\left(x-y\right) -\mu}{\sigma}\right)H_\lambda\left(x,y\right)dxdy\\
    &+ \iint_{\left(-b_m^2,\infty\right)^2} x^{k_1}y^{k_2}e^{-ax-by}\phi\left(\frac{\left(x-y\right) -\mu}{\sigma}\right)H_\lambda\left(x,y\right)\frac{1}{\Phi\left(\min\left(x,y\right)\right)^2}dxdy.
\end{align*}
The first integral is finite by Lemma 1. To see that the second integral is finite, we split the integral into the subsets $\left\{x\leq y\right\}$ and $\left\{y<x\right\}$. Over the region including $\{x < 0, x \leq y\}$, we have the bound $\Phi\left(\min\left(x,y\right)\right) = 1-\Phi(-x) \geq \phi(x)\frac{x^2-1}{-x^3}$. With this bound, we have

\begin{align*}
    &\iint_{\mathbb{R}^2\setminus \left(-b_m^2,\infty\right)^2 \cap \left\{x \leq y\right\}} x^{k_1}y^{k_2}e^{-ax-by}\phi\left(\frac{\left(x-y\right) -\mu}{\sigma}\right)H_\lambda\left(x,y\right)\frac{1}{\Phi\left(\min\left(x,y\right)\right)^2}dxdy\\
    \leq &\iint_{\mathbb{R}^2\setminus \left(-b_m^2,\infty\right)^2 \cap \left\{x \leq y\right\}} x^{k_1}y^{k_2}e^{-ax-by}\phi\left(\frac{\left(x-y\right) -\mu}{\sigma}\right)H_\lambda\left(x,y\right)\frac{1}{\left[\phi\left(x\right)\frac{x^2-1}{-x^3}\right]^2}dxdy\\
    \leq &2\pi\iint_{\mathbb{R}^2\setminus \left(-b_m^2,\infty\right)^2 \cap \left\{x \leq y\right\}} x^{k_1}y^{k_2}e^{-ax-by+x^2}\phi\left(\frac{\left(x-y\right) -\mu}{\sigma}\right)e^{-\frac{1}{2}e^{-x}}x^2dxdy\\
    \leq&2\pi\sum_{i=0}^{k_2}\binom{k_2}{i}\int_{-\infty}^{-b_m^2}e^{-\left(a+b\right)x+x^2}e^{-\frac{1}{2}e^{-x}}x^{2+k_1+k_2-i}dx\int^{\infty}_{-\infty}\left(-\delta\right)^ie^{b\delta}\phi\left(\frac{\delta - \mu}{\sigma}\right)d\delta,
\end{align*}
where $\delta = x-y$. Finally, to see the same holds for $\left\{y<x\right\}$, note that $x$ and $y$ are symmetric, and redefine $\delta = y-x$.

\qed

\begin{lemma}\label{appendix2}
    Integrals of the form 
    \begin{multline}
        \iint_{\mathbb{R}^2}x^{k_1}y^{k_2}e^{-ax-by}\phi\left(\textcolor{black}{q_m(y,x)}\right)e^{-\frac{x^2}{2b_m^2}}\Phi_{\rho_m}\left(u_m\left(x\right),u_m\left(y\right)\right)^{m-2}\\\left(m-1\right)\sum_{i=2}^{\infty}\frac{\left(1-\Phi_{\rho_m}\left(u_m\left(x\right),u_m\left(y\right)\right)\right)^i}{i}dxdy
    \end{multline} are finite for all sufficiently large $m$.
\end{lemma}

\noindent \textit{Proof.}
Assume $m$ large enough such that $\rho_m>0$. Using the bound on the infinite sum from the previous lemma (and once again dropping the $(e^{-x}+e^{-y})^2$ term), we re-write the integral as

\begin{align*}
    \iint_{\mathbb{R}^2}x^{k_1}y^{k_2}e^{-ax-by}\phi\left(\textcolor{black}{q_m(y,x)}\right)e^{-\frac{x^2}{2b_m^2}}\Phi_{\rho_m}\left(u_m\left(x\right),u_m\left(y\right)\right)^{m-3}dxdy
\end{align*}
The only difference between the integral above and the integral in \textcolor{black}{Lemma}~\ref{Convergence of expectations} is the difference between the power $m-2$ and $m-3$, which makes little material difference, as well as the inclusion of the powers of $x,y$. As such the proof that this integral is finite is similar to the proof in \textcolor{black}{Lemma}~\ref{Convergence of expectations}. Define $l, j, C > 0$ such that $|\delta|^{i} < Ce^{l\delta}+Ce^{-j\delta}$ for all $i \in \left\{0,...,k_2\right\}$. Using this bound, we get

\begin{align*}
    &\iint_{\mathbb{R}^2}x^{k_1}y^{k_2}e^{-ax-by}\phi\left(\textcolor{black}{q_m(y,x)}\right)e^{-\frac{x^2}{2b_m^2}}\Phi_{\rho_m}\left(u_m\left(x\right),u_m\left(y\right)\right)^{m-3}dxdy \\
    &\leq \int_{\mathbb{R}}x^{k_1}e^{-\left(a+b\right)x}e^{-\frac{x^2}{2b_m^2}}\Phi\left(u_m\left(x\right)\right)^{m-3}\int_{\mathbb{R}}\left(x-\delta\right)^{k_2}e^{k_2\delta}\phi\left(\frac{\delta - \left(1-\rho_m\right)\left(b_m^2-x\right)}{b_m\sqrt{1-\rho_m^2}}\right)d\delta dx\\
    &\leq C\sum_{i=0}^{k_2}\binom{k_2}{i}\int_{\mathbb{R}}x^{k_1+i}e^{-\left(a+b\right)x-\frac{x^2}{2b_m^2}}\Phi\left(u_m\left(x\right)\right)^{m-3}\int_{\mathbb{R}}\left(e^{l\delta}+e^{-j\delta}\right)e^{k_2\delta}\phi\left(\frac{\delta - \left(1-\rho_m\right)\left(b_m^2-x\right)}{b_m\sqrt{1-\rho_m^2}}\right)d\delta dx\\
    &= b_m\sqrt{1-\rho_m^2}C\sum_{i=0}^{k_2}\binom{k_2}{i}\int_{\mathbb{R}}x^{k_1+i}e^{-\left(a+b\right)x-\frac{x^2}{2b_m^2}}\Phi\left(u_m\left(x\right)\right)^{m-3}\\
    &\hspace{15mm}\left[e^{\left(1-\rho_m\right)\left(b_m^2-x\right)\left(l+k_2\right) + \frac{\left(l+k_2\right)^2}{2}\left(b_m^2\left(1-\rho_m^2\right)\right)} + e^{\left(1-\rho_m\right)\left(b_m^2-x\right)\left(k_2-j\right) + \frac{\left(k_2-j\right)^2}{2}\left(b_m^2\left(1-\rho_m^2\right)\right)}\right]dx.
\end{align*}

Again, we note that $\rho_m \rightarrow 1$ and $b_m\sqrt{1-\rho_m^2} \rightarrow 2\lambda$. Thus for large enough $m$ such that $\left(1-\rho_m\right)\left(k_2+l\right) -\left(a+b\right) < 0$ and $\left(1-\rho_m\right)\left(k_2-j\right) -\left(a+b\right) < 0$, we write the remaining integral as

\begin{align*}
    \int_{\mathbb{R}}x^{k_1+i}e^{-tx-\frac{x^2}{2b_m^2}}\Phi_{\rho_m}\left(u_m\left(x\right),u_m\left(y\right)\right)^{m-3}dx
\end{align*}
for some $t > 0$. Over the interval $\left[0,\infty\right)$

\begin{align*}
    &\int_0^\infty x^{k_1+i}e^{-tx-\frac{x^2}{2b_m^2}}\Phi_{\rho_m}\left(u_m\left(x\right),u_m\left(y\right)\right)^{m-3}dx
    \leq\int_0^\infty  x^{k_1+i}e^{-tx}dx = \frac{\Gamma\left(k_1+i+1\right)}{t^{k_1+i-1}}.
\end{align*}

Over the interval $\left(\left.-\infty,0\right]\right.$, there exist $C_2, v > 0$ such that $x^{k_2+i} < Ce^{-vx}$. The integral becomes

\begin{align*}
    &\int_{-\infty}^0 Ce^{-\left(t+v\right)x-\frac{x^2}{2b_m^2}}\Phi_{\rho_m}\left(u_m\left(x\right),u_m\left(y\right)\right)^{m-3}dx.
\end{align*}

The proof that this term is finite is similar to that presented in \textcolor{black}{Lemma}~\ref{Convergence of expectations}, with the power $m-2$ instead of $m-3$. Thus this term is finite.

\qed

\begin{lemma}\label{appendix3}
    Integrals of the form 
    \begin{equation*}
        \iint_{\mathbb{R}}x^{k_1}y^{k_2}e^{-ax-by}\gamma\phi\left(\gamma\right)H_\lambda\left(x,y\right)dxdy
    \end{equation*}
    where $\gamma$ is defined to be the second order Lagrange remainder of the Taylor series expansion of $\Phi\left(q_m(x,y)\right)$ about $\lambda + \frac{y-x}{2\lambda}$, or the same term defined symmetrically about x and y, are finite.
\end{lemma}

\noindent \textit{Proof.}

Without loss of generality, we prove the lemma for $\gamma$ defined as the Lagrange remainder of the Taylor series expansion of $\Phi\left(q_m(x,y)\right)$ about $\lambda + \frac{y-x}{2\lambda}$, as the proof with the alternate definition is similar. Recall that \textcolor{black}{
\begin{equation*}
    \min\left(\textcolor{black}{q_m(y,x)},\lambda + \frac{y-x}{2\lambda}\right) \leq \gamma \leq \max\left(\textcolor{black}{q_m(y,x)},\lambda + \frac{y-x}{2\lambda}\right).
\end{equation*}}Thus, we know that $|\gamma|$ is always bounded by some polynomial in $x$ and $y$. Continuing, we absorb this term into the terms $x^{k_1}$ and $y^{k_2}$. If $\textcolor{black}{q_m(y,x)}$ and $\lambda + \frac{y-x}{2\lambda}$ have the same sign, then we may write either \textcolor{black}{ 
\begin{equation*}
    \phi\left(\gamma\right) < \phi\left(\lambda - \frac{y-x}{2\lambda}\right) \quad \text{or} \quad \phi\left(\gamma\right) < \phi\left(\textcolor{black}{q_m(y,x)}\right).
\end{equation*}}In either case, we see that the integral over this region is finite by the same argument used in Lemma~\ref{appendix1}. These terms have the same sign when $\delta := x-y > \max\left(2\lambda^2, b_m^2\left(1-\rho_m\right) + x\sqrt{1-\rho_m^2}\right)$ or $\delta < \min\left(2\lambda^2, b_m^2\left(1-\rho_m\right) + x\sqrt{1-\rho_m^2}\right)$. \textcolor{black}{Define $R_m$ to be the region in $(x,y)$ when $\delta$ is between these bounds.} We have the result:

\begin{align*}
    &\iint_{\textcolor{black}{R_m}}x^{k_1}y^{k_2}e^{-ax-by}\gamma\phi\left(\gamma\right)H_\lambda\left(x,y\right)dxdy\\
    &\leq\int x^{k_1}e^{-\left(a+b\right)x}e^{-\frac{1}{2}e^{-x}}\int^{\max\left(2\lambda^2, b_m^2\left(1-\rho_m\right) + x\sqrt{1-\rho_m^2}\right)}_{\min\left(2\lambda^2, b_m^2\left(1-\rho_m\right) + x\sqrt{1-\rho_m^2}\right)}\left(x-\delta\right)^{k_2}e^{b\delta}dxd\delta\\
    &\leq C\sum_{i=0}^{k_2}\binom{k_2}{i}\int x^{k_1+i}e^{-\left(a+b\right)x}e^{-\frac{1}{2}e^{-x}}\int^{\max\left(2\lambda^2, b_m^2\left(1-\rho_m\right) + x\sqrt{1-\rho_m^2}\right)}_{\min\left(2\lambda^2, b_m^2\left(1-\rho_m\right) + x\sqrt{1-\rho_m^2}\right)}e^{\left(b+j\right)\delta}+e^{\left(b-l\right)\delta}dxd\delta\\
    &\leq C\sum_{i=0}^{k_2}\binom{k_2}{i}\int x^{k_1+i}e^{-\left(a+b\right)x}e^{-\frac{1}{2}e^{-x}}\left(e^{\left(b+j\right)\max\left(2\lambda^2, b_m^2\left(1-\rho_m\right) + x\sqrt{1-\rho_m^2}\right)}+e^{\left(b-l\right)\max\left(2\lambda^2, b_m^2\left(1-\rho_m\right) + x\sqrt{1-\rho_m^2}\right)}\right.\\
    &\left.\hspace{5mm}-e^{\left(b+j\right)\min\left(2\lambda^2, b_m^2\left(1-\rho_m\right) + x\sqrt{1-\rho_m^2}\right)}-e^{\left(b-l\right)\min\left(2\lambda^2, b_m^2\left(1-\rho_m\right) + x\sqrt{1-\rho_m^2}\right)}\right)dx
\end{align*}
for some sufficiently large $j,l$. This integral is surely finite for $m$ large enough such that $\left(b+j\right)x\sqrt{1-\rho_m^2} < \left(a+b\right)$ by the same argument as in Lemma~\ref{integral}.

\qed

\begin{lemma}\label{appendix4}
    Integrals of the form
    \begin{equation*}
        \iint e^{-ax-by}x^{k_1}y^{k_2}\Phi\left(\lambda + \frac{x-y}{2\lambda}\right)\Phi\left(\lambda + \frac{y-x}{2\lambda}\right)H_{\lambda}\left(x,y\right)dxdy
    \end{equation*}
    or 
    \begin{equation*}
        \iint e^{-ax-by}x^{k_1}y^{k_2}\Phi\left(\lambda + \frac{x-y}{2\lambda}\right)\Phi\left(\lambda + \frac{y-x}{2\lambda}\right)e^{-\frac{x^2}{2b_m^2}}\Phi_{\rho_m}\left(u_m\left(x\right),u_m\left(y\right)\right)^{m-2}dxdy
    \end{equation*}
    are finite for $a+b \geq 1$, $\textcolor{black}{k_1,k_2 \in \mathbb{N}}$.
\end{lemma}

\noindent\textit{Proof.}
The key to this proof is the bound
\begin{equation*}
   \frac{\phi\left(x\right)}{x}\left(1-\frac{1}{x^2}\right) < 1-\Phi\left(x\right) <  \frac{\phi\left(x\right)}{x} 
\end{equation*}
for $x>0$ as stated in 26.2.12 in \textcites{abramowitz+stegun}. Rewriting the first equation as

\begin{align*}
    &\iint e^{-ax-by}x^{k_1}y^{k_2}\Phi\left(\lambda + \frac{x-y}{2\lambda}\right)\Phi\left(\lambda + \frac{y-x}{2\lambda}\right)H_{\lambda}\left(x,y\right)dxdy\\
    &\leq \sum_{i=0}^{k_2}\binom{k_2}{i}\int e^{-\left(a+b\right)x}x^{k_1+i}e^{-\frac{1}{2}e^{-x}}dx\int e^{b\delta}\delta^{k_2-i}\Phi\left(\lambda + \frac{\delta}{2\lambda}\right)\Phi\left(\lambda - \frac{\delta}{2\lambda}\right)d\delta\\
     &= \sum_{i=0}^{k_2}\binom{k_2}{i}\int e^{-\left(a+b\right)x}x^{k_1+i}e^{-\frac{1}{2}e^{-x}}dx\int e^{b\delta}\delta^{k_2-i}\left(1-\Phi\left(-\lambda - \frac{\delta}{2\lambda}\right)\right)\left(1-\Phi\left(-\lambda + \frac{\delta}{2\lambda}\right)\right)d\delta\\
     &= \sum_{i=0}^{k_2}\binom{k_2}{i}\int e^{-\left(a+b\right)x}x^{k_1+i}e^{-\frac{1}{2}e^{-x}}dx\left[\int^{-2\lambda\left(\lambda+1\right)}_{-\infty} e^{b\delta}\delta^{k_2-i}\phi\left(\lambda + \frac{\delta}{2\lambda}\right)d\delta\right.\\
     &\left.\hspace{10mm}+ \int_{2\lambda\left(\lambda+1\right)}^{\infty} e^{b\delta}\delta^{k_2-i}\phi\left(\lambda - \frac{\delta}{2\lambda}\right)d\delta + \int_{-2\lambda\left(\lambda+1\right)}^{2\lambda\left(\lambda+1\right)} e^{b\delta}\delta^{k_2-i}d\delta  \right].
\end{align*}

The rest of the proof for this form of the integral is similar to Lemma~\ref{integral}. The proof for the second form is a combination of this argument and that used in Lemma~\ref{appendix2}.

We also note that inclusion or exclusion of the term $\left(m-1\right)\sum_{i=2}^{\infty}\frac{\left(1-\Phi_{\rho_m}\left(u_m\left(x\right),u_m\left(y\right)\right)\right)^i}{i}$ does not affect the integrability, which can be shown by the same technique used in Lemma~\ref{appendix1}.

\begin{redlemma}
    All terms in (\ref{Split Expectation 1}) that were excluded from the body of the proof of \textcolor{black}{Lemma}~\ref{bias-variance theorem} necessarily converge to zero faster than $\frac{1}{\ln\left(m\right)}$ or $\lambda_m - \lambda$, and are integrable.
\end{redlemma}
We begin with the left over terms from (\ref{Split Expectation 1}). The first terms we handle are

\begin{enumerate}[I]
    \item $\displaystyle=\frac{1}{4b_m^4}\int_y^{\infty}\Phi\left(\lambda + \frac{x-z}{2\lambda}\right)e^{-z-\nu_z}z^4dz$
    \item $\displaystyle=\int_{y}^{\infty}\phi\left(\lambda + \frac{x-z}{2\lambda}\right)\left[\textcolor{black}{q_m(x,z)}-\lambda -\frac{x-z}{2\lambda}\right]e^{-z}\left[-\frac{z^2}{2b_m^2} + e^{-\nu_z}\frac{z^4}{4b_m^4}\right]dz$
    \item $\displaystyle=\int_{y}^{\infty}\gamma_z\phi\left(\gamma_z\right)\left[\textcolor{black}{q_m(x,z)}-\lambda -\frac{x-z}{2\lambda}\right]^2e^{-z}\left[1-\frac{z^2}{2b_m^2} + e^{-\nu_z}\frac{z^4}{4b_m^4}\right]dz$
    \item $\displaystyle=\frac{1}{4b_m^4}\int_x^{\infty}e^{-t-\nu_t}t^4dt$
    \item $\displaystyle=\frac{1}{m}e^{-x}$
    \item $\displaystyle=\frac{1}{m}\int_y^{\infty}\Phi\left(\lambda + \frac{x-z}{2\lambda}\right)e^{-z}dz$
    \item $\displaystyle=-\left(m-1\right)\sum_{i=2}^{\infty}\frac{\left(1-\Phi_{\rho_m}\left(u_m\left(x\right),u_m\left(y\right)\right)\right)^i}{i}$
\end{enumerate}
Our goal for terms I to VI is to bound each term by the form $\frac{1}{g\left(m\right)}e^{-ax-by}x^{k_1}y^{k_2}$ for some $a,b,k_1,k_2 \in \mathbb{N} \cup \left\{0\right\}$ and $g$ such that $\frac{1}{g\left(m\right)}$ is dominated by at least one of $\frac{1}{b_m^2}$ or $\lambda_n - \lambda$.  Recalling the definition of $\nu_z$ and $\nu_t$ to be the second order Lagrange remainder term for the expansion of $e^{-\frac{z^2}{2b_m^2}}$ or $e^{-\frac{t^2}{2b_m^2}}$ about 1, we know $0 < \nu_t, \nu_z$ and thus $e^{-\nu_t}, e^{-\nu_z} < 1$. Thus we drop the term wherever it may appear. We see that
\begin{align*}
    \text{I} &=\frac{1}{4b_m^4}\int_y^{\infty}\Phi\left(\lambda + \frac{x-z}{2\lambda}\right)e^{-z-\nu_z}z^4dz\\
    &\leq \frac{1}{4b_m^4}\int_y^{\infty}e^{-z}z^4dz\\
    & = e^{-y}\frac{y^4+4y^3+12y^2+24y+24}{4b_m^4}\\
    \text{II} & =\int_{y}^{\infty}\phi\left(\lambda + \frac{x-z}{2\lambda}\right)\left[\textcolor{black}{q_m(x,z)}-\lambda -\frac{x-z}{2\lambda}\right]e^{-z}\left[-\frac{z^2}{2b_m^2} + e^{-\nu_z}\frac{z^4}{4b_m^4}\right]dz\\
    &\leq \int_{y}^{\infty}\left[\lambda_n - \lambda + \left(x-z\right)\left(\frac{1}{b_m\sqrt{1-\rho_m^2}} - \frac{1}{2\lambda}\right) + z \sqrt{\frac{1-\rho_m}{1+\rho_m}}b_m^{-1}\right]e^{-z}\left[\frac{z^2}{2b_m^2} + \frac{z^4}{4b_m^4}\right]dz\\
    &=e^{-y}\left(\lambda_n - \lambda + x\left(\frac{1}{b_m\sqrt{1-\rho_m^2}} - \frac{1}{2\lambda}\right)\left(\frac{y^2+2y+2}{2b_m^2} + \frac{y^4+4y^3+12y^2+24y+24}{4b_m^4}\right)\right)\\
    & \hspace{10mm}- e^{-y}\left[\left(\frac{1}{b_m\sqrt{1-\rho_m^2}} - \frac{1}{2\lambda}\right) + \sqrt{\frac{1-\rho_m}{1+\rho_m}}b_m^{-1}\right]\left( + \frac{y^3+3y^2+6y+6}{2b_m^2}\right.\\
    &\hspace{15mm}+\left. \frac{y^5+5y^4+20y^3+60y^2+120y+120}{4b_m^4}\right)\\
    \text{IV} &= \frac{1}{4b_m^4}\int_x^{\infty}e^{-t-\nu_t}t^4dt \\
    &\leq \frac{1}{4b_m^4}e^{-x}\left(x^4+4x^3+12x^2+24x+24\right)\\
    \text{VI} &= \frac{1}{m}\int_y^{\infty}\Phi\left(\lambda + 
    \frac{x-z}{2\lambda}\right)e^{-z}dz\\
    &\leq \frac{1}{m}e^{-y}
\end{align*}
We note that term V is already in the form required. To handle III, recall the definition of $\gamma_z$ as the Lagrange remainder from the expansion of $\Phi\left(\textcolor{black}{q_m(x,z)}\right)$ about $\lambda + \frac{x-z}{2\lambda}$. Then we know 
\textcolor{black}{\begin{equation*}
    |\gamma_z| < \left|\max\left\{\textcolor{black}{q_m(x,z)},\lambda + \frac{x-z}{2\lambda}\right\}\right|.
\end{equation*}}Thus, including this term only includes another polynomial in $x$ and $y$. Using this polynomial bound, and dropping $\phi\left(\gamma_z\right) < 1$, we see that the analysis of III is similar to that of II. 
Thus, since every term is of the form $\frac{1}{g\left(m\right)}e^{-ax-by}x^{k_1}y^{k_2}$, the integral will take the form
\begin{equation}
\iint_{\mathbb{R}^2}\frac{\partial}{\partial\lambda}\ln\left(h_\lambda\left(x,y\right)\right)H_\lambda\left(x,y\right)\phi\left(\lambda + \frac{y-x}{2\lambda}\right)\frac{1}{2\lambda}\frac{1}{g\left(m\right)}e^{-ax-by}x^{k_1}y^{k_2}dxdy.
\end{equation}
Recall that each term in $\frac{\partial}{\partial\lambda}\ln\left(h_\lambda\left(x,y\right)\right)$ is also of the form $e^{-ax-by}x^{k_1}y^{k_2}$. Thus each of these integrals is finite and converges to zero by Lemma~\ref{integral}. Of the terms listed above, only VII remains, which is explicitly taken care of in Lemma~\ref{appendix1}. There yet remains one set of term we need to handle:
\begin{multline*}
    e^{\zeta_a}\left(f\left(x,y,m\right) - \frac{1}{2}\textcolor{black}{q_m(y,x)}^2
    + \frac{1}{2}\left(\lambda + \frac{y-x}{2\lambda}\right)^2 - \ln\left(\frac{b_m\sqrt{1-\rho_m^2}}{2\lambda}\right)\right)^2.
\end{multline*}

We note that this is the product of any two terms (or square of any term) analyzed either here or in the body of \textcolor{black}{Lemma}~\ref{bias-variance theorem} multiplied by the term $e^{\zeta_a}$. Thus, so long as the integral stays finite, each term will necessarily converge to zero faster than either $\frac{1}{\ln\left(m\right)}$ or $\lambda_n - \lambda$. Recall the definition of $\zeta_a$ as the second order Lagrange remainder term for the expansion of 
\begin{equation*}
    \exp\left\{f\left(x,y,m\right)-\frac{x^2}{2b_m^2} - \frac{1}{2}\textcolor{black}{q_m(y,x)}^2 + \frac{1}{2}\left(\lambda + \frac{y-x}{2\lambda}\right)^2 - \ln\left(\frac{b_m\sqrt{1-\rho_m^2}}{2\lambda}\right)\right\}
\end{equation*}
about 0. If $\zeta_a < 0$, then we may simply drop the term and the analysis of each term is just as before. However, if 
\textcolor{black}{\begin{equation*}
    0 > \zeta_a > f\left(x,y,m\right)-\frac{x^2}{2b_m^2} - \frac{1}{2}\textcolor{black}{q_m(y,x)}^2 + \frac{1}{2}\left(\lambda + \frac{y-x}{2\lambda}\right)^2 - \ln\left(\frac{b_m\sqrt{1-\rho_m^2}}{2\lambda}\right)
\end{equation*}}then our integrand includes the term
\begin{align*}
    &H_\lambda\left(x,y\right)e^{-x}\phi\left(\lambda + \frac{y-x}{2\lambda}\right)\frac{1}{2\lambda}\exp\left\{f\left(x,y,m\right)-\frac{x^2}{2b_m^2} - \frac{1}{2}\textcolor{black}{q_m(y,x)}^2 + \frac{1}{2}\left(\lambda + \frac{y-x}{2\lambda}\right)^2\right.\\
    & \hspace{10mm}\left.-\ln\left(\frac{b_m\sqrt{1-\rho_m^2}}{2\lambda}\right)\right\}\\
    &=\Phi_{\rho_m}\left(u_m\left(x\right),u_m\left(y\right)\right)^{m-1}e^{-x}\phi\left(\textcolor{black}{q_m(y,x)}\right)\frac{1}{b_m\sqrt{1-\rho_m^2}}.
\end{align*}

As such, our integrals will take the form
\begin{equation*}
    \iint_{\mathbb{R}^2}e^{-ax-by}x^{k_1}y^{k_2}e^{-\frac{x^2}{2b_m^2}}\Phi_{\rho_m}\left(u_m\left(x\right),u_m\left(y\right)\right)^{m-1}\phi\left(\textcolor{black}{q_m(y,x)}\right)\frac{1}{b_m\sqrt{1-\rho_m^2}}dxdy.
\end{equation*}
or
\begin{multline*}
    \iint_{\mathbb{R}^2}e^{-ax-by}x^{k_1}y^{k_2}e^{-\frac{x^2}{2b_m^2}}\Phi_{\rho_m}\left(u_m\left(x\right),u_m\left(y\right)\right)^{m-1}\phi\left(\textcolor{black}{q_m(y,x)}\right)\frac{1}{b_m\sqrt{1-\rho_m^2}}\\\left(m-1\right)\sum_{i=2}^{\infty}\frac{\left(1-\Phi_{\rho_m}\left(u_m\left(x\right),u_m\left(y\right)\right)\right)^i}{i}dxdy\textcolor{black}{.}
\end{multline*}
Both of these cases are handled by Lemma~\ref{appendix2}. Thus, all terms converge to 0 faster than $\frac{1}{\ln\left(m\right)}$ or $\lambda_n -\lambda$.

\qed

\begin{redlemma}
    All terms in (\ref{Split Expectation 2}) that were excluded from the body of the proof of \textcolor{black}{Lemma}~\ref{bias-variance theorem} necessarily converge to zero faster than $\frac{1}{ln\left(m\right)}$ or $\lambda_m - \lambda$, and are integrable.
\end{redlemma}

\noindent\textit{Proof.}
The first set of excluded terms we consider are:

\begin{enumerate}[I]
    \item $=\phi\left(\lambda + \frac{y-x}{2\lambda}\right)\left(\textcolor{black}{q_m(y,x)} - \lambda - \frac{y-x}{2\lambda}\right)\phi\left(\lambda + \frac{x-y}{2\lambda}\right)\left(\textcolor{black}{q_m(x,y)} - \lambda - \frac{x-y}{2\lambda}\right)$
    \item $=\left(\textcolor{black}{q_m(x,y)} - \lambda - \frac{x-y}{2\lambda}\right)^2\gamma_2\phi\left(\gamma_2\right)\Phi\left(\lambda + \frac{y-x}{2\lambda}\right)$
    \item $=\left(\textcolor{black}{q_m(x,y)} - \lambda - \frac{x-y}{2\lambda}\right)^2\gamma_2\phi\left(\gamma_2\right)\phi\left(\lambda + \frac{y-x}{2\lambda}\right)\left(\textcolor{black}{q_m(y,x)} - \lambda - \frac{y-x}{2\lambda}\right)$
    \item $=\left(\textcolor{black}{q_m(y,x)} - \lambda - \frac{y-x}{2\lambda}\right)^2\gamma_1\phi\left(\gamma_1\right)\phi\left(\lambda + \frac{y-x}{2\lambda}\right)\left(\textcolor{black}{q_m(y,x)} - \lambda - \frac{y-x}{2\lambda}\right)$
    \item $=\left(\textcolor{black}{q_m(y,x)} - \lambda - \frac{y-x}{2\lambda}\right)^2\gamma_1\phi\left(\gamma_1\right)\Phi\left(\lambda + \frac{x-y}{2\lambda}\right)$
    \item $=\left(\textcolor{black}{q_m(y,x)} - \lambda - \frac{y-x}{2\lambda}\right)^2\gamma_1\phi\left(\gamma_1\right)\left(\textcolor{black}{q_m(x,y)} - \lambda - \frac{x-y}{2\lambda}\right)^2\gamma_2\phi\left(\gamma_2\right)$
\end{enumerate}

Recall that we have shown that $\textcolor{black}{q_m(x,y)} - \lambda - \frac{x-y}{2\lambda}$ and $\textcolor{black}{q_m(y,x)} - \lambda - \frac{y-x}{2\lambda}$ are polynomials in $x$ and $y$. Additionally, we know that $|\gamma_1|,|\gamma_2|$ are bounded by a polynomial in $x$,$y$. We can drop any $\Phi$ above as they must be between 0 and 1. Thus terms I, III, and IV all contain a term $\phi\left(\lambda + \frac{y-x}{2\lambda}\right)$ or $\phi\left(\lambda + \frac{x-y}{2\lambda}\right)$ and are thus integrable by Lemma~\ref{integral}. Terms II,V and VI all depend on either $\phi\left(\gamma_1\right)$ or $\phi\left(\gamma_2\right)$ and are thus integrable by Lemma~\ref{appendix3}. We note that each of these terms must converge faster than $\frac{1}{ln\left(m\right)}$ or $\lambda_n - \lambda$ as each term includes a product of $\textcolor{black}{q_m(x,y)} - \lambda - \frac{x-y}{2\lambda}$ and/or $\textcolor{black}{q_m(y,x)} - \lambda - \frac{y-x}{2\lambda}$. As shown, each term in these polynomials must converge at least as fast as one of $\frac{1}{ln\left(m\right)}$ or $\lambda_n - \lambda$ and thus their product must converge faster.

The next bundle of terms is the product of any term in $f'\left(x,y,m\right)$ with any one of the 6 terms above or

\begin{enumerate}[(a)]
    \item $=\Phi\left(\lambda + \frac{x-y}{2\lambda}\right)\Phi\left(\lambda + \frac{y-x}{2\lambda}\right)$
    \item $=\Phi\left(\lambda + \frac{x-y}{2\lambda}\right)\phi\left(\lambda + \frac{y-x}{2\lambda}\right)\left(\textcolor{black}{q_m(y,x)} - \lambda - \frac{y-x}{2\lambda}\right)$
    \item $=\Phi\left(\lambda + \frac{y-x}{2\lambda}\right)\phi\left(\lambda + \frac{x-y}{2\lambda}\right)\left(\textcolor{black}{q_m(x,y)} - \lambda - \frac{x-y}{2\lambda}\right)$
\end{enumerate}

Recall that all terms in $f'\left(x,y,m\right)$ are of the form $e^{-ax-by}x^{-k_1}y^{-k_2}$. Thus, whichever Lemma showed that the terms by themselves were finite also hold when multiplied by a term in $f'\left(x,y,m\right)$. Any term including the infinite sum $\left(m-2\right)\sum_{i=2}^{\infty}\frac{\left(1-\Phi_{\rho_m}\left(u_m\left(x\right),u_m\left(y\right)\right)\right)^i}{i}$ is also finite by the analogous lemma in the appendix.

Finally, we address the term including $e^{\zeta_b}$. Just as with $\zeta_a$, if $\zeta_b < 0$, we may simply drop the term, and all terms are finite by the same analysis above. If, however, $\zeta_b$ is positive, we have the bound

\begin{equation*}
    e^{\zeta_b} < e^{-\frac{x^2}{2b_m^2}-\frac{y^2}{2b_m^2}}\frac{\Phi_{\rho_m}\left(u_m\left(x\right),u_m\left(y\right)\right)^{m-2}}{H_\lambda\left(x,y\right)}.
\end{equation*}
Thus our integrals become

\begin{equation*}
    \iint_{\mathbb{R}^2}\frac{\partial}{\partial \lambda}\ln\left(h_\lambda\left(x,y\right)\right)e^{-\frac{x^2}{2b_m^2}-\frac{y^2}{2b_m^2}}\Phi_{\rho_m}\left(u_m\left(x\right),u_m\left(y\right)\right)^{m-2}Xdxdy,
\end{equation*}
where $X$ is any included term from above. Again, any term in $f\left(x,y,m\right)$ form $e^{-ax-by}x^{k_1}y^{k_2}$, so the integrability of each term only depend on which of the 9 terms is included in the product. Each term will be finite by one of the lemmas given. \qed

Finally, we give the form of the bias term $A$. \textcolor{black}{Recall the definitions}
\begin{equation*}
L_1 := \lim_{m\rightarrow\infty}\frac{\sqrt{k}}{b_m^2}
     \quad\text{and}\quad
    L_2:= \lim_{m\rightarrow\infty}\sqrt{k}(\lambda-\lambda_m)
\end{equation*}
and assume both are finite. $A$ is given by:

\begin{align*}
    &L_1\iint_{\mathbb{R}^2}\frac{\partial}{\partial\lambda}\ln(h_\lambda(x,y))H_\lambda(x,y)\phi\left(\lambda + \frac{y-x}{2\lambda}\right)\frac{1}{2\lambda}e^{-x}\left(.5 \int_y^\infty\Phi\left(\lambda + \frac{x-z}{2\lambda}\right)e^{-z}z^2dz\right.\\
    &\hspace{5mm}\left.+.5e^{-x}(x^2+2x+2) -.5x^2 + \lambda^2 - \frac{\lambda}{2}\int_y^\infty\phi\left(\lambda + \frac{x-z}{2\lambda}\right)e^{-z}(x+z)dz\right)dxdy\\
    +&L_2\iint_{\mathbb{R}^2}\frac{\partial}{\partial\lambda}\ln(h_\lambda(x,y))H_\lambda(x,y)\phi\left(\lambda + \frac{y-x}{2\lambda}\right)\frac{1}{2\lambda}e^{-x}\left(\frac{(x-y)^2}{4\lambda^3}-\lambda+\frac{1}{\lambda}\right.\\
    &\hspace{5mm}-\left.e^{-x}\int_y^\infty\phi\left(\lambda+\frac{x-z}{2\lambda}\right)\left(\frac{x-z}{2\lambda^2}-1\right)dz\right)dxdy\\
    +&L_1\iint_{\mathbb{R}^2}\frac{\partial}{\partial\lambda}\ln(h_\lambda(x,y))H_\lambda(x,y)e^{-x-y}\Phi\left(\lambda + \frac{x-y}{2\lambda}\right)\Phi\left(\lambda + \frac{y-x}{2\lambda}\right)\left. .5 \int_y^\infty\Phi\left(\lambda + \frac{x-z}{2\lambda}\right)e^{-x}z^2dz\right.\\
    &\hspace{5mm}-\frac{\lambda}{2}\int_y^\infty\phi\left(\lambda + \frac{x-z}{2\lambda}\right)e^{-z}(x+z)dz + .5e^{-x}(x^2+2x+2) -x^2\\
    &\left.\hspace{5mm}+\lambda\Phi\left(\lambda + \frac{y-x}{2\lambda}\right)^{-1}\phi\left(\lambda + \frac{y-x}{2\lambda}\right)(x+y)\right)dxdy \\
    +&L_2\iint_{\mathbb{R}^2}\frac{\partial}{\partial\lambda}\ln(h_\lambda(x,y))H_\lambda(x,y)e^{-x-y}\Phi\left(\lambda + \frac{x-y}{2\lambda}\right)\left(\phi(\lambda + \frac{y-x}{2\lambda})\left(\frac{y-x}{\lambda^2}-2\right)\right.\\
    &\left.\hspace{5mm}-\Phi\left(\lambda + \frac{y-x}{2\lambda}\right)e^{-x}\int_y^\infty\phi\left(\lambda+\frac{x-z}{2\lambda}\right)\left(\frac{x-z}{2\lambda^2}-1\right)dz\right)dxdy.
\end{align*}

\end{document}